\documentclass[reqno]{amsart}
\usepackage{amssymb, stmaryrd}
\usepackage{amscd}
\usepackage{pdfsync} 

\newtheorem{theorem}{Theorem}
\newtheorem{lemma}[theorem]{Lemma}

\newtheorem{example}[theorem]{Example}
\newtheorem{remark}[theorem]{Remark}
\newtheorem{proposition}[theorem]{Proposition}
\newtheorem{corollary}[theorem]{Corollary}

\newcommand{\F}{\boldsymbol{F}}
\newcommand{\FF}{\overline{\boldsymbol{F}}}
\newcommand{\TT}{\boldsymbol{\mathcal{T}} }
\newcommand{\UU}{\boldsymbol{\mathcal{U}} }
\newcommand{\MM}{\boldsymbol{\mathcal{M}} }
\newcommand{\WW}{\boldsymbol{\mathcal{W}} }
\newcommand{\VV}{\boldsymbol{\mathcal{V}} }
\newcommand{\DD}{\boldsymbol{\mathcal{D}} }
\newcommand{\KK}{\boldsymbol{\mathcal{K}} }

\newcommand{\Fc}{\mathcal{F}}

\usepackage{color}  
\definecolor{darkgreen}{rgb}{0.03, 0.5, 0.03}
 \newcommand{\ec}{\color{black}} %

 \newcommand{\co} {\boldsymbol c}
 \newcommand{\Na} {\mbox{\rm Na}}
  \newcommand{\Kr} {\mbox{\rm Kr}}
  \newcommand{\Max} {\mbox{\rm Max}}
    \newcommand{\Spec} {\mbox{\rm Spec}}
      \newcommand{\QMax} {\mbox{\rm QMax}}
        \newcommand{\QSpec} {\mbox{\rm QSpec}}

\newcommand{\cal} {\mathcal}

\newcommand{\Fb}{\boldsymbol{\overline{F}}}
\newcommand{\f}{\boldsymbol{{f}}}

  \newcommand{\stf} {\star{_{\!{_f}}}}
    \newcommand{\stt} {\widetilde{\star}}
      \newcommand{\bo}{{\boldsymbol{\circledast}}_0}
       \newcommand{\ba}{{ \boldsymbol{\circledast}}}

       \newcommand{\bto}{\boldsymbol{(\wedge_{\TT})}}
       \newcommand{\bag}{\boldsymbol{\langle\wedge_{\TT}\rangle}}
       \newcommand{\bsq}{ \boldsymbol{[\wedge_{\TT}]}}

          \newcommand{\btow}{\boldsymbol{(\wedge_{\WW})}}
 \newcommand{\btov}{\boldsymbol{(\wedge_{\VV})}}
  \newcommand{\btod}{\boldsymbol{(\wedge_{\DD})}}
   \newcommand{\btok}{\boldsymbol{(\wedge_{\KK})}}

           \newcommand{\bagv}{\boldsymbol{\langle\wedge_{\VV}\rangle}}
                 \newcommand{\bagd}{\boldsymbol{\langle\wedge_{\DD}\rangle}}
         \newcommand{\bagk}{\boldsymbol{\langle\wedge_{\KK}\rangle}}

        \newcommand{\bsqv}{ \boldsymbol{[\wedge_{\VV}]}}
            \newcommand{\bsqd}{ \boldsymbol{[\wedge_{\DD}]}}
        \newcommand{\bsqk}{ \boldsymbol{[\wedge_{\KK}]}}

          \newcommand{\btom}{\boldsymbol{(\wedge_{\!\MM})}}
       \newcommand{\bagm}{\boldsymbol{\langle\wedge_{\!\MM}\rangle}}
       \newcommand{\bsqm}{ \boldsymbol{[\wedge_{\!\MM}]}}

\begin{document}

\title[Polynomial extensions of star and semistar operations]
{An overring-theoretic approach to polynomial extensions of star and semistar operations}
{\author[G. W. Chang and M. Fontana]{Gyu Whan Chang and Marco Fontana}}

\date{\today}

\address{Dipartimento di Matematica, Universit\`a degli Studi
``Roma Tre'', 00146 Rome, Italy.}
\email{fontana@mat.uniroma3.it }
\address{Department of Mathematics, University of Incheon, Incheon 406- 840, Korea}
\email{whan@incheon.ac.kr}

\begin{abstract}
Call a semistar operation $\ba$ on   the  polynomial domain $D[X]$ an extension (respectively,
a strict extension) of a semistar operation $\star$ defined on  an  integral domain $D$, with quotient field $K$, if $E^\star = (E[X])^{\ba}\cap K$
(respectively, $E^\star [X]= (E[X])^{\ba}$) for all nonzero $D$-submodules $E$ of $K$. In this paper,
we study the general properties of the above defined extensions and link our work with earlier efforts, centered on the stable semistar operation case,
at defining semistar operations on $D[X]$ that are ``canonical'' extensions (or, ``canonical'' strict extensions) of  semistar
operations on $D$.
 \end{abstract}

\maketitle

\bigskip
\section{Background results}

Let $D$ be an integral domain with quotient field $K$. Let
$\boldsymbol{\overline{F}}(D)$ denote the set of all nonzero
$D$--submodules of $K$ and let $\boldsymbol{F}(D)$ be the set of
all nonzero fractional ideals of $D$, i.e., $E \in
\boldsymbol{F}(D)$ if $E \in \boldsymbol{\overline{F}}(D)$ and
there exists a nonzero $d \in D$ with $dE \subseteq D$. Let
$\f(D)$ be the set of all nonzero finitely generated
$D$--submodules of $K$. Then, obviously $\f(D)
\subseteq \boldsymbol{F}(D) \subseteq
\boldsymbol{\overline{F}}(D)$.

Following Okabe-Matsuda \cite{o-m}, a \emph{semistar operation} on $D$ is a map $\star:
\boldsymbol{\overline{F}}(D) \to \boldsymbol{\overline{F}}(D),\ E \mapsto E^\star$,  such that, for all $x \in K$, $x \neq 0$, and
for all $E,F \in \boldsymbol{\overline{F}}(D)$, the following
properties hold:
\begin{enumerate}
\item[$(\star_1)$] $(xE)^\star=xE^\star$;
 \item[$(\star_2)$] $E
\subseteq F$ implies $E^\star \subseteq F^\star$;
\item[$(\star_3)$] $E \subseteq E^\star$ and $E^{\star \star} :=
\left(E^\star \right)^\star=E^\star$.
\end{enumerate}

 The semistar operation defined by $E^\star = K$ for all $E \in \boldsymbol{\overline{F}}(D)$ is called the {\it trivial semistar operation on $D$} and it is denoted by $e_D$ (or, simply, by $e$).
A \emph{(semi)star operation} is a semistar operation that,
restricted to $\boldsymbol{F}(D)$,  is a star operation (in the
sense of \cite[Section 32]{G}). It is easy to see that a
semistar operation $\star$ on $D$ is a (semi)star operation if and
only if $D^\star = D$.

 If $\star$ is a semistar operation on $D$, then we can
consider a map\ $\star_{\!_f}: \boldsymbol{\overline{F}}(D) \to
\boldsymbol{\overline{F}}(D)$ defined, for each $E \in
\boldsymbol{\overline{F}}(D)$, as follows:

\centerline{$E^{\star_{\!_f}}:=\bigcup \{F^\star\mid \ F \in
\boldsymbol{f}(D) \mbox{ and } F \subseteq E\}$.}

\noindent It is easy to see that $\star_{\!_f}$ is a semistar
operation on $D$, called \emph{the semistar operation of finite
type associated to $\star$}.  Note that, for each $F \in
\boldsymbol{f}(D)$, $F^\star=F^{\star_{\!_f}}$.  A semistar
operation $\star$ is called a \emph{semistar operation of finite
type} if $\star=\star_{\!_f}$.  It is easy to see that
$(\star_{\!_f}\!)_{\!_f}=\star_{\!_f}$ (that is, $\star_{\!_f}$ is
of finite type).

If $\star_1$ and $\star_2$ are two semistar operations on $D$, we
say that $\star_1 \leq \star_2$ if $E^{\star_1} \subseteq
E^{\star_2}$, for each $E \in \Fb(D)$. This is equivalent to say
that $\left(E^{\star_{1}}\right)^{\star_{2}} = E^{\star_2}=
\left(E^{\star_{2}}\right)^{\star_{1}}$, for each $E \in
\Fb(D)$.  Obviously, for each semistar operation  $\star$,
we
have $\star_{\!_f} \leq \star$. Let $d_D$ (or, simply, $d$)  be the \it identity (semi)star operation on $D$,   \rm clearly $d \leq \star$, for all semistar operation $\star$ on $D$.

We say that a nonzero ideal $I$ of $D$ is a
\emph{quasi-$\star$-ideal} if $I^\star \cap D = I$, a
\emph{quasi-$\star$-prime  ideal} if it is a prime quasi-$\star$-ideal,
and a \emph{quasi-$\star$-maximal  ideal} if it is maximal in the set of
all   proper   quasi-$\star$-ideals. A quasi-$\star$-maximal ideal is  a
prime ideal. It is possible  to prove that each  proper   quasi-$\star_{_{\!
f}}$-ideal is contained in a quasi-$\star_{_{\! f}}$-maximal
ideal.  More details can be found in \cite[page 4781]{FL}. We
will denote by $\QMax^{\star}(D)$  (respectively, $\QSpec^\star(D)$) the set of the
quasi-$\star$-maximal ideals  (respectively, quasi-$\star$-prime ideals) of $D$.
When $\star$ is a (semi)star operation, the notion of  quasi-$\star$-ideal coincides with the ``classical'' notion of  \it  $\star$-ideal \rm (i.e., a nonzero ideal $I$ such that $I^\star = I$).

\smallskip

  If $\Delta$ is a set of prime ideals of an integral domain
  $D$,  then the semistar operation $\star_\Delta$ defined
  on $D$ as follows

  \centerline{$
  E^{\star_\Delta} := \bigcap \{ED_P \;|\;\, P \in \Delta\}\,,
  \;  \textrm {  for each}    \; E \in \boldsymbol{\overline{F}}(D)\,,
  $}

  \noindent is called \it the spectral semistar operation   on $D$ associated to
  \rm
  $\Delta$.
  A semistar operation $\star$ on an integral domain $D$ is
  called
  \it a
  spectral semistar operation \rm if there exists a subset $
  \Delta$ of the prime spectrum of $D$, $\mbox{\rm Spec}(D)$,  such that $\,\star =
  \star_\Delta\,$.

  When $\Delta := \QMax^{\star_{_{\! f}}}(D)$, we set $\stt:= \star_{\Delta}$, i.e.

  \centerline{
  $E^{\stt} := \bigcap \left \{ED_P \mid P \in  \QMax^{\star_{_{\! f}}}(D) \right\}$,  \;  for each $E \in \boldsymbol{\overline{F}}(D)$.
  }

 A semistar operation $\star$ is \emph{stable} if $(E \cap F)^\star
= E^\star \cap F^\star$, for each $E,F \in \Fb(D)$.
Spectral semistar operations are stable \cite[Lemma 4.1(3)]{FH}. In particular, $\stt$ is a semistar operation stable and of finite type; and,  conversely, if a semistar operation $\star$ is stable and of finite type then $\star = \stt$ \cite[Corollary 3.9(2)]{FH}.

\smallskip
\smallskip
Let $T$ be an overring of an integral domain $D$, let $\iota: D \hookrightarrow T$ be the canonical embedding
and  let $\star$ be a semistar operation on $D$. We denote by  $\star_\iota$   the semistar operation on $T$ defined by $E^{\star_\iota} := E^\star$,  for each $E \in \Fb(T) \ (\subseteq \Fb(D))$.   \\
  Conversely, let $\ast$ be a semistar operation on $T$ and let $\ast^\iota$ be the semistar operation on $D$ defined by $E^{\ast^\iota} := (ET)^\ast$,  for each $E \in \Fb(D)$.\\
   It is not difficult to see that
${(\ast^{\iota})}{_{_{\! f}}} = ({\ast}{_{_{\!f}}})^{\iota}$  and  if $\star$ is a semistar operation of finite type (respectively, a stable semistar operation) on $D$ then ${\star_\iota}$ is a semistar operation of finite type (respectively, a stable semistar operation) on $T$ (cf. for instance \cite[Proposition 2.8]{FL1} and \cite[Propositions 2.11 and 2.13]{pi}). \\
Clearly, if $\star= d_D$ then $(d_D)_\iota = d_T$. In case $\ast = d_T$, the semistar operation of finite type $(d_T)^\iota$ (defined by $E \mapsto ET$ for all $E \in \Fb(D)$) is denoted also by $\star_{\{T\}}$ and it is stable if and only if $T$ is a flat overring of $D$ \cite[Proposition 1.7]{uda} and \cite [Theorem 7.4(1)]{matsumura}.

\smallskip

 By $v_D$ (or, simply, by $v$) we denote  the $v$--(semi)star
operation defined as usual by  $E^v := (D:(D:E))$, for each $E\in
\boldsymbol{\overline{F}}(D)$. By  $t_D$ (or, simply, by $t$) we
denote  $(v_D)_{_{\! f}}$ the $t$--(semi)star operation on $D$ and
by  $w_D$ (or just by $w$) the stable semistar operation of finite
type associated to $v_D$ (or, equivalently, to $t_D$), considered
by F.G.  Wang  and R.L. McCasland in \cite{WMc97} (cf. also \cite{gv});  i.e.,  $w_D :=
\widetilde{\ v_D} = \widetilde{\ t_D}$.  Clearly $w_D\leq t_D \leq v_D$.  Moreover, it is easy to see that for each   (semi)star operation $\star$ on $D$, we have $\star \leq v_D$ and $\stf \leq t_D$ (cf. also \cite[Theorem 34.1(4)]{G}).

\smallskip

We recall from \cite[Chapter V]{FHP97} (see also \cite[Chapter 4]{Po}) that \it a localizing
system of ideals of $D$ \rm  is a family $\Fc$ of ideals of $D$
such that:

\begin{enumerate}
\item[{\bf(LS1)}] If  $I \in \Fc$ and $J$ is an ideal of $D$ such
that $I \subseteq J$, then $J \in \Fc$.
  \item[{\bf(LS2)}] If $I \in \Fc$ and $J$ is an ideal of $D$ such that
  $(J:_DiD) \in  \Fc$, for each $i \in I$, then $J \in \Fc$.
\end{enumerate}

A localizing system $\Fc$ is \emph{finitely generated} if, for each
$I \in \Fc$, there exists a finitely generated ideal $J \in \Fc$
such that $J \subseteq I$.

The relation between  stable semistar operations and localizing
systems has been deeply investigated by M. Fontana and J. Huckaba
in \cite{FH} and by F. Halter-Koch in the context of module
systems  \cite{HK01}.    In the following proposition, we summarize some of the results that we need  (see \cite[Proposition~2.8,
Proposition 3.2, Proposition 2.4, Corollary 2.11, Theorem 2.10
(B)]{FH}).

\begin{proposition}
 \label{prop:loc1}  Let $D$ be an integral domain.
\begin{enumerate}

\item[\bf (1)\rm] If $\star$ is a semistar operation on $D$, then
$\Fc^\star:=\left\{I \mbox{ ideal of $D$ } \mid I^\star = D^\star
\right\}$ is a localizing system (called \emph{the localizing
system associated to $\star$}).
 \item[\bf (2)\rm] If $\star$ is a semistar operation of finite type, then
$\Fc^\star$ is a finitely generated localizing system.
  \item[\bf (3)\rm] Let $\star_\Fc$ or, simply, ${\overline{\star}}$ be the semistar operation associated to a given  localizing system $\Fc$ of $D$ and defined by $E \mapsto E^{\overline{\star}}:=
\bigcup \left \{(E:J) \,\vert\, J\in\Fc \right \}$, for each $E \in
\Fb(D)$.  Then  $\star_\Fc$ (called \emph{the semistar operation
associated to the localizing system $ \Fc$}) is  a stable semistar
operation on $D$.
  \item[\bf (4)\rm] ${\overline{\star}} \leq \star$ and $\Fc^{\star} = \Fc^{\overline{\star}}$.
    \item[\bf (5)\rm] $\overline{\star} = \star$ if and only if $\star$ is
  stable.
  \item[\bf (6)\rm] If $\Fc$ is a finitely
generated localizing system, then $\star_\Fc$ is a finite type
(stable) semistar operation.

  \item[\bf (7)\rm]  $\Fc^{\star_{_{\! f}}}= (\Fc^\star){_{_{\! f}}} := \{ I \in \Fc^\star \mid I\supseteq J, \mbox{ for some finitely generated ideal } J\in \Fc^\star\} $ and  $\widetilde{\star}=\overline{\ \!\star_{_{\! f}}}$, i.e., $ \widetilde{\star}$ is the stable semistar operation of finite type associated to the localizing
system $\Fc^{\star_{_{\! f}}}$. In particular, for each $E \in \Fb(D)$, we   have:
$$
E^{\stt} = \bigcup \{(E:J) \mid J \in \f(D),\  J \subseteq D,  \mbox{ and } J^\star =D^\star \}\,.
$$
   \item[\bf (8)\rm] If $\Fc^\prime$ and $\Fc^{\prime \prime}$ are
two localizing systems of $D$, then ${\Fc}^\prime \subseteq
\Fc^{\prime \prime}$ if and only if $\star_{_{{\Fc}^\prime}} \leq
\star_{_{{\Fc}^{\prime \prime}}}$. \hfill $\Box$
\end{enumerate}
\end{proposition}

If $I$ is a nonzero fractional ideal of $D$, we say that $I$ is
\emph{$\star$--invertible} if $(II^{-1})^\star = D^\star$.   From the
definitions and from the fact that $\QMax^{\stf}(D) = \QMax^{\stt}(D)$ \cite[Corollary   3.5(2)]{FL} it follows easily   that a nonzero fractional ideal $I$ is
$\tilde{\star}$--invertible if and only if $I$ is $\star_{_{\!
f}}$--invertible.  An integral domain $D$ is called a \emph{Pr\"ufer $\star$--multiplication domain} (for short, \emph{P$\star$MD}) if each $I \in \f(D)$ is $\star_{_{\! f}}$--invertible. It is easy to see that the notions of P$v$MD, P$t$MD and P$w$MD coincide. Obviously, a P$d$MD is a Pr\"ufer domain, and conversely \cite[Theorem 22.1]{G}.

\smallskip

If $R$ is a ring (not necessarily an integral domain) and $X$ an indeterminate over $R$, then the ring
$R(X):=\{f/g \, \vert \; f,g \in R[X] \mbox{ and }
\boldsymbol{c}(g)=R \}$ (where $\boldsymbol{c}(g)$ is the content
of the polynomial $g$) is called the \emph{Nagata ring} of $R$
\cite[Proposition 33.1]{G}.

In case of an integral domain equipped with a semistar operation, we have a general ``semistar version'' of the Nagata ring. The following result was proved in \cite[Propositions 3.1 and 3.4]{FL}
(cf.  also \cite[Proposition 2.1]{K-89}).

\begin{proposition}\label{nagata}
Let $\star$ be a nontrivial semistar operation on an integral
domain $D$. Set ${\mathcal N}^\star := {\mathcal N}_D^\star:=\{h \in D[X] \mid h
\neq 0 \mbox{ and } \boldsymbol{c}(h)^\star=D^\star \}$ and
$$
\Na(D, \star) := D[X]_{{\mathcal N}^\star}\,.
$$
 Then,
\begin{enumerate}
\item[{\bf (1)}] ${\mathcal N}^\star $ is a saturated multiplicative subset of
$D[X]$ and ${\mathcal N}^\star ={\mathcal N}^{\star_{\!_f}}=D[X] \smallsetminus \bigcup \{Q[X]
\, \vert \; Q \in \QMax^{\star_{\!_f}}(D) \}$.
\item[{\bf (2)}] $\Max(\Na(D, \star) )= \{Q[X]_{{\mathcal N}^\star } \mid  Q \in
 \QMax^{\star_{\!_f}}(D) \}$ and $ \QMax^{\star_{\!_f}}(D)$ coincides with
the canonical image in\ $\Spec(D)$ of\
$\Max \left(\Na(D, \star) \right)$.
\item[{\bf (3)}] $\Na(D, \star) =
\bigcap \{D_Q(X) \, \vert \; Q \in  \QMax^{\star_{\!_f}}(D) \}$.
   \item[{\bf (4)}] For each $E \in \Fb(D)$, $E^{\stt} = E\Na(D, \star) \cap K$. \hfill $\Box$
\end{enumerate}
\end{proposition}

Let
     $\star$ be a semistar operation on $D$. If $F$ is in $\boldsymbol{f}(D)$, we say that
     $F$ is
            \it  $\star$--\texttt{eab} \rm (respectively,  \it
             $\star$--\texttt{ab}\rm)
  if
$(FG)^{\star}
             \subseteq (FH)^{\star}$ implies that $G^{\star}\subseteq H^{\star} $, with $G,\ H \in
\boldsymbol{f}(D)$, (respectively,  with $G,\ H \in
\overline {\boldsymbol{F}}(D)$).

      An operation $\star$ is said to be   \it \texttt{eab} \rm (respectively, \it \texttt{ab}\rm\ \!)  if each $F\in \boldsymbol{f}(D)$ is $\star$--\texttt{eab}  (respectively, $\star$--\texttt{ab}). \  An \texttt{ab} operation is obviously an \texttt{eab} operation.
      We  note  that if $\star$ is an \texttt{eab} semistar operation then $\, {\star_{_{\! f}} }\,$ is also an \texttt{eab} semistar operation, since they agree on all finitely generated ideals.
 Let $\star$ be a semistar operation of finite type, then
$ \star$  is an \texttt{eab} semistar operation if and only if  $\star$ is an \texttt{ab}  semistar operation. In this situation, we say that $\star$ is an \emph{\texttt{(e)ab} semistar operation}. In particular, from the previous result it follows that the notions of $\star$--\texttt{eab} semistar operation and   $\star_{_{\! f}}$--\texttt{(e)ab}  semistar operation coincide \cite[Lemma 3 and Proposition 4]{FL-09}.

Given an arbitrary semistar operation $\star$ on an integral domain
$D$, it is possible to associate to $\star$, an \texttt{eab} semistar
operation of finite type $\, \star_a \, $ of $D$, called \it the \texttt{ab}
semistar operation associated to $\star$, \rm defined as follows:

\centerline{$ F^{\star_a} :=
 \cup\{((FH)^\star:H)\,\; |\; \, H \in \boldsymbol{f}(D) \}, \; \; \textrm {for each } \,
 F \in \boldsymbol{f}(D) \, ,
$}

\noindent and,  in general,

\centerline{$ E^{\star_a} := \cup\{F^{\star_a} \,\; |\; \, F \subseteq E\,,\; F
\in \boldsymbol{f}(D) \},\; \; \textrm {for each } \, E \in
\boldsymbol{\overline{F}}(D), $}

\noindent \cite [Definition 4.4]{FL1}.  Note that if $\star$ is an
\texttt{(e)ab} semistar operation of finite type then $\star = \star_{a}$, and conversely
\cite [Proposition 4.5]{FL1}.
  More information on the operation ${\star_a} $, introduced for ideal systems in \cite[page 41]{Jaffard} (see also Lorenzen's original paper \cite{Lz}) can be found in  \cite{OM1},
\cite{o-m}, \cite{Koch:1997}, \cite{HK98}, and \cite{FL0}.

Let $\star$ be a semistar operation on $D$ and
let $V$ be a valuation overring of $D$.  We say that $V$ is a \it $\star$--valuation overring of
$D$ \rm if, for each $F \in \boldsymbol{f}(D)\,$, $ F^\star \subseteq FV\,$ (or equivalently,
$\star_f \leq \star_{\{V\}}$.
Note that a valuation overring $V$ of $D$ is a
$\star$--valuation overring of $D$ if and only if $V^{\star_{f}} =V$.
 More details of semistar valuation overrings can be found in
\cite{FL0}, \cite{FL1} (cf.  also \cite{Jaffard},  \cite{Koch:1997} and
\cite{HK03}).

  \begin{proposition}\label{prop:Kr} \rm
  \cite[Proposition 3.3, Theorem 3.11, Theorem 5.1, Corollary 5.2, Corollary
  5.3]{FL1}, \cite[Theorem 3.5]{FL0}.
  \it Let $\star$ be any semistar
  operation defined on an integral domain $D$ with quotient field $K$ and
  let $\star_{a}$ be the \texttt{ab}  semistar operation associated to $\star$.   Set
$$
\begin {array} {rl}
\Kr(D,\star) := \{ f/g \; \,|\, & f,g \in D[X] \setminus \{0\} \;\;
\mbox{ \it and there exists } \; h \in D[X] \setminus \{0\} \; \\ &
\mbox{ \it such that } \; (\boldsymbol{c}(f)\boldsymbol{c}(h))^\star
\subseteq (\boldsymbol{c}(g)\boldsymbol{c}(h))^\star \,\} \, \cup\,
\{0\}\,.
\end{array}
  $$
  Then, we have:
\begin{enumerate}
\item[{\bf (1)}] $\Kr(D,\star)$ is a B\'ezout domain with quotient field
$K(X)\,,$ called \rm the Kronecker function ring of $D$ with respect to
the semistar operation $\star\,.$ \it
\item[{\bf (2)}]  ${\Na}(D,\star) \subseteq {\Kr}(D,\star)\,.$
\item[{\bf (3)}]  ${\Kr}(D,\star) = {\Kr}(D,\star_a)\,$.
\item[{\bf (4)}]
$ E^{\star_{a}} = E{\Kr}(D,\star) \cap K \,,$ \; for each $ E \in
\boldsymbol{\overline{F}}(D)$\,.
\item[{\bf (5)}]\ $\Kr(D, \star) = \bigcap \{ V(X) \mid V \mbox{ is a $\
\star$--valuation overring of } D \}\,.$
\item[{\bf (6)}] If $\,F := (a_{0},a_{1},\ldots, a_{n}) \in
 \boldsymbol{f}(D)$\, and $\,f(X) :=a_{0}+ a_{1}X +\ldots +a_{n}X^n \in
 K[X]\,,$ then:

 \centerline{ $ F{\Kr}(D,\star) = f(X){\Kr}(D,\star) =
 {\boldsymbol{c}}(f){\Kr}(D,\star)\,.  $}

 \end{enumerate} \vskip -0.5cm \hfill $\square$

  \end{proposition} \rm

When $\star =d$, the $d$--valuation overrings of $D$ are just the valuation overrings of $D$. In this case, we set:
$$ \Kr(D) := \Kr(D, d) =  \bigcap \{ V(X) \mid V \mbox{ is a
valuation overring of } D \}\,.$$
Moreover, if we   denote by $b_D$  (or, simply, by $b$) the \texttt{ab} semistar operation of finite type $(d_D)_a$ then, for each $E \in \Fb(D)$,
$$
E^b = E\Kr(D) \cap K = \bigcap \{EV \mid  V \mbox{ is a valuation overring of } D \}\,.
$$

\begin{remark} \label{prufer}
\rm Recall  that a Pr\"ufer domain $D$ can be characterized by the fact that each $F \in \f(D)$ is invertible. Since an invertible ideal is always a $v$--ideal (and, in particular, a $t$--ideal), then the following are equivalent \cite[Theorem 24.7 and Theorem 34.1(4)]{G}: \it
\begin{enumerate}

\item[{\bf (i)}]  $D$ is a Pr\"ufer domain;

\item[{\bf (ii)}]  $D$ is integrally closed and $d = b$;

\item[{\bf (iii)}]  $D$ is integrally closed and $d = t$.

\end{enumerate}\rm
\end{remark}

\section{Results}


 Let $D$ be an integral domain with quotient field $K$, let $X$ be an indeterminate over $K$. We start with some basic facts. Note that some of the statements contained in the following result were also proved in \cite[Proposition 2.1]{M} and, in the star operation setting, in \cite[Propositions 2.1 and 2.2]{HMM}.

\begin{lemma} \label{ast-zero}
Given a semistar operation $\ba$  on $D[X]$,  for each $E \in \FF(D)$ set:
$$ E^{\bo} := (E[X])^{\ba} \cap K\,.$$
Then:
\begin{enumerate}
\item[{\bf (1)}] $ \bo$ is a semistar operation on $D$ called \rm the semistar operation canonically induced by $\ba$ on $D$. \it  In particular, if $\ba$ is a (semi)star operation on $D[X]$, then $\bo$ is a (semi)star operation on $D$.
\item[{\bf (2)}] \it  $(E^{\bo}[X])^{\ba} = (E[X])^{\ba}$ for all $E \in \FF(D)$.
\item[{\bf (3)}] $({\ba}_{_{\! f}})_0 = (\bo)_{_{\! f}}$.

\item[{\bf (4)}] \it If $\ba$ is a semistar operation of finite type (respectively, stable), then $\bo$ is a   semistar  operation of finite type (respectively, stable).

\item[{\bf (5)}] \it   If $\ba^\prime$ and $\ba^{\prime\prime}$ are two semistar operations on $D[X]$ and $\ba^\prime \leq \ba^{\prime\prime}$, then \ $\ba^\prime_0 \leq \ba^{\prime\prime}_0$.

 \item[{\bf (6)}]
   $(\widetilde{\ba})_0 = \widetilde{\ \bo\ }$.

\item[{\bf (7)}] \it $(d_{D[X]})_{_{\!{0}}} = d_{D}$, \ $(w_{D[X]})_{_{\!{0}}} = w_{D}$, \ $(t_{D[X]})_{_{\!{0}}} = t_{D}$, \ $(v_{D[X]})_{_{\!{0}}} = v_{D}$, \ and  $(b_{D[X]})_{_{\!{0}}} = b_{D}$.

\end{enumerate} \rm
\end{lemma}
\begin{proof} {\bf (1)} Set $\star:=\bo$.
 It is easy to see that, if $E \in
\FF(D)$, then $E\subseteq
E^{\star}$ and if $ E_{1}, E_{2} \in
\FF(D)$ with  $E_{1}\subseteq E_{2}$, then $ E_{1}^{\star}
\subseteq  E_{2}^{\star}$. Moreover:
$$
\begin{array}{rl}
\left(E^{{\star}}\right)^{{\star}} & \hskip -5pt =
(((E[X])^{\ba} \cap K)[X])^{\ba} \cap K \\
 & \hskip -5pt  \subseteq  ((E[X])^{\ba}[X] \cap K[X])^{\ba} \cap K\\
& \hskip -5pt =((E[X])^{\ba} \cap K[X])^{\ba}
\cap K\\
& \hskip -5pt \subseteq \left((E[X])^{\ba}\right)^{\ba}
\cap K = (E[X])^{\ba} \cap K  =
E^{{\star}} \,.
 \end{array}
 $$

Thus $\left(E^{{\star}}\right)^{{\star}} = E^{{\star}}$.
\ec Moreover, for each nonzero $z \in K$, we have:
 $$
\begin{array}{rl}
zE^{{\star}}& \hskip -5pt =
z((E[X])^{\ba} \cap K) \\
 & \hskip -5pt =(z(E[X])^{\ba} \cap zK) =(z(E[X])^{\ba} \cap K)\\
& \hskip -5pt =(zE[X])^{\ba} \cap K=
\left(zE\right)^{{{\star}}} \,.
 \end{array}
$$

\noindent  In particular, if $(D[X])^{\ba} = D[X]$, then $D^{\star} = (D[X])^{\ba} \cap K
= D[X] \cap K = D$.

{\bf (2)}  Note that $E[X] \subseteq E^{\bo}[X]= ((E[X])^{\ba} \cap K)[X]   \subseteq (E[X])^{\ba} \cap K[X] \subseteq (E[X])^{\ba} $. Therefore, $(E[X])^{\ba} \subseteq (E^{\bo}[X])^{\ba}\subseteq ((E[X])^{\ba})^{\ba} =  (E[X])^{\ba}$.

 {\bf (3)}  Let $z \in E^{({\ba}_{_{\! f}})_0}$.
 Then there exists $F\in  \f(D[X])$ such that $F \subseteq E[X]$
  and $z \in F^{\ba} \cap K$. Let $I :=\co_D(F)$.
  Clearly,  $I \in \f(D)$, $I \subseteq E$ and $F\subseteq I[X] \subseteq E[X]$. Therefore, $z \in F^{\ba} \cap K \subseteq (I[X])^{\ba} \cap K = I^{\bo}$, and so $z \in E^{({\bo})_{\! f}}$.  Conversely, if $z \in E^{({\bo})_{\! f}}$, then $z \in I^{\bo} = (I[X])^{\ba} \cap K $ for some  $I \in \f(D)$, $I \subseteq E$. This implies that $z \in E^{({\ba}_{_{\! f}})_0}$.

  {\bf (4)} The ``finite type part'' is a particular case of (3). The ``stable part" is a straightforward consequence of the definitions.

   {\bf (5)}  is straightforward.

    {\bf (6)}     Clearly,
     $(\widetilde{\ba})_0 \leq \widetilde{\ \bo\ }$, since  $(\widetilde{\ba})_0 \leq {\ \bo\ }$ by (5) and $(\widetilde{\ba})_0 $ is a stable operation of finite type by  statement (4). Let  $a \in E^{\widetilde{\ \bo \ }}$ with $E \in \FF(D)$. Then, there exists a nonzero finitely generated ideal $J$ of $D$ such that $J^{\bo} = D^{\bo}$ and $aJ \subseteq E$.  On the other hand,
      $(J[X])^{\ba}=  (J^{\bo}[X])^{\ba} =  (D^{\bo}[X])^{\ba} =(D[X])^{\ba}$ by (2).
     Since $J[X]$ is a nonzero finitely  generated ideal of $D[X]$, $(J[X])^{\ba} = (D[X])^{\ba}$ and $aJ[X] \subseteq E[X]$, then  $a \in (E[X])^{\widetilde{\ba}}$ and so $a \in (E[X])^{\widetilde{\ba}} \cap K = E^{({\widetilde{\ba}})_0}$.

{\bf (7)}  The statement for the $d$--operations is trivial. For the $w$--, $t$--, and $v$--operation,   it
 is an easy consequence of the following equalities  \cite[Proposition 4.3]{HH} for all fractional ideals $E\in \F(D)$:
  $$
  (E[X])^{v_{D[X]}} = E^{v_D}[X], \;\;    (E[X])^{t_{D[X]}} = E^{t_D}[X],  \;\;    (E[X])^{w_{D[X]}} = E^{w_D}[X]\,.
  $$
  (Note that if $E\in \FF(D)\setminus \F(D)$, then  $E[X] \in  \FF(D[X])\setminus \F(D[X])$, and so $E^{v_D}[X] = K[X]$ and $(E[X])^{v_{D[X]}}= K(X)$, however  $ (E[X])^{v_{D[X]}} \cap K = K = E^{v_D}$.)

 We want to prove next that
 $$
 (E[X])^{b_{D[X]}} = E^{b_D}[X], \mbox{  for all } E\in \FF(D)\,.
 $$We use the fact that $E^{b_D} = \bigcup \{ (EI:I) \mid I\in \f(D) \}$ (respectively,
 $(E[X])^{b_{D[X]}} = \bigcup \{ (E[X]F:F) \mid F\in \f(D[X]) \}$) (see \cite[page 349]{ZS} and \cite[Section 19.3]{HK}).  Let $z \in E^{b_D}[X] \ (\subseteq K[X])$. Then  $z \in (EI:I)[X] = (E[X]I[X]:I[X])$ for some $I \in \f(D)$, and so, in particular, $z \in   \bigcup \{ (E[X]F:F) \mid F\in \f(D[X]) \} $.  Conversely, let $z \in   \bigcup \{ (E[X]F:F) \mid F\in \f(D[X]) \} \ (\subseteq K[X])$. Then $zF \subseteq E[X]F$, and so (by \cite[Theorem 28.1]{G})
 %
 %
 there exists a positive integer $m$
 such that $\co_D(z)(\co_D(z)^m\co_D(F)) = \co_D(z)^{m+1}\co_D(F)$ =
 $\co_D(z)^m\co_D(zF) \subseteq \co_D(z)^m\co_D(E[X]F) \subseteq \co_D(z)^m\co_D(E[X])\co_D(F) $ $=E(\co_D(z)^m\co_D(F))$,
where  $\co_D(z)^m\co_D(F) \in \f(D)$. Therefore, $\co_D(z) \subseteq E^{b_D}$ and thus $z \in \co_D(z)[X] \subseteq E^{b_D}[X]$.
 \end{proof}

  \begin{remark} \rm Note that  the equality $(E[X])^{b_{D[X]}} = E^{b_D}[X], \mbox{  for all } E\in \FF(D)$ proved in (7) of  Lemma \ref{ast-zero}, is equivalent to each of the following equalities:
    $$
    \begin{array}{rl}
  & E[X]\Kr(D[X], b_{D[X]}) \cap K(X)  = (E\Kr(D, b_{D}) \cap K)[X] = E\Kr(D, b_{D}) \cap K[X]\,, \\
  & \bigcap \{E[X] W \mid W \mbox{ valuation overring of } D[X]\}  = \\
   &  \hskip 160pt    \left(\bigcap \{E V \mid V \mbox{ valuation overring of } D\}\right)\![X]\,.
 \end{array}
  $$
   \end{remark}

\begin{remark} \label{CF}
  \rm Given an arbitrary  multiplicative subset $\mathcal S$   of $D[X]$, Chang and Fontana   in \cite{CF}  investigated  the map $E \mapsto ED[X]_{\mathcal{S}}\cap K$,\  defined for all $E \in \overline{\boldsymbol{F}}(D)$, showing that it gives rise to a semistar operation $\star$  on  $D$, having the properties that $D^\star = R :=D[X]_{\mathcal{S}}\cap K$, and that $R$ is $t$-linked to $(D, \star)$ (i.e.,  for each nonzero finitely generated ideal $I$ of $D$, $I^{\star}=D^\star$ implies $(IR)^{t_R}= R$  \cite[Section 3]{eBF};
or, equivalently,
$R = R^{\stt}$ \cite[Lemma 2.9]{CF2}.)    One of the main results of the paper by Chang and Fontana \cite[Theorem 2.1]{CF} is recalled below, since it is strictly linked to the theme of the present work.

Note that, to a   multiplicative subset $\mathcal S$   of $D[X]$, we can associate the semistar operation $\ba_{\mathcal S}$ on $D[X]$ defined by $A^{\ba_{\mathcal S}} := A_{\mathcal S} = \bigcup \{(A: J)\mid J \mbox{ ideal of } D[X],\ J \cap {\mathcal S} \neq \emptyset \}=  AD[X]_{\mathcal{S}}$,
for all $A \in \FF(D[X])$ \cite[Proposition 2.10]{FH}.  Therefore, by  Lemma \ref{ast-zero}, we obtain immediately that the map $E \mapsto ED[X]_{\mathcal{S}}\cap K=: E^{\circlearrowleft_{\mathcal{S}}} $,\  defined for all $E \in \FF(D)$, gives rise to a semistar operation ${\circlearrowleft_{\mathcal{S}}}$ on $D$ coinciding with $({\ba_{\mathcal S}})_0$.
Clearly, if $\overline{\mathcal S} :=  D[X] \setminus  \bigcup \{ Q \mid  Q \in \Spec(D[X]) \mbox{ and } Q \cap \mathcal S = \emptyset \} $  is the saturation of the multiplicative set $\mathcal S$, then ${\ba_{\mathcal S}}  = {\ba_{\overline{\mathcal S}}} $ and so, in particular, ${\circlearrowleft_{\mathcal{S}}} = {\circlearrowleft_{\overline{\mathcal S}}}$.

In order to deepen our knowledge of the semistar operation ${\circlearrowleft_{\mathcal{S}}}$, we
need a definition of a stronger version of saturation. Set:
$$\mathcal S^{\sharp} := D[X] \setminus  \bigcup \{ P[X] \mid  P \in \Spec(D) \mbox{ and }  P[X] \cap \mathcal S = \emptyset \}.
$$
  It is clear that $\mathcal S^{\sharp}$ is a saturated multiplicative set of $D[X]$ and that  $\mathcal S^{\sharp}$ contains the saturation of $\mathcal S$, i.e., $ \mathcal S^{\sharp} \supseteq  \overline{\mathcal S} \supseteq {\mathcal S} $. We  call $\mathcal S^{\sharp}$ \it the extended saturation of ${\mathcal{S}}$ in $D[X]$ \rm and a multiplicative set  $\mathcal S$ of $D[X]$ is called \it extended saturated \rm  if $\mathcal S=\mathcal S^{\sharp}$.    Clearly, in general,  ${\ba_{\mathcal S^{\sharp}}}  \geq  {\ba_{{\mathcal S}}} \ (  = {\ba_{\overline{\mathcal S}}}) $. However, it can be shown  that  $({\ba_{\mathcal S^{\sharp}}} )_0 =  ({\ba_{{\mathcal S}}} )_0$.  For this, let
  $$
  \boldsymbol{\Delta}:= \boldsymbol{\Delta}(\mathcal{S}) := \{ P \in \Spec(D)\mid    P[X] \cap \mathcal{S} = \emptyset \}\,;
  $$
  obviously,  $\boldsymbol{\Delta}(\mathcal{S}) = \boldsymbol{\Delta}(\mathcal{S}^{\sharp})$. Let  $\boldsymbol{\nabla}:= \boldsymbol{\nabla}(\mathcal{S})$ be the set of the maximal elements of $\boldsymbol{\Delta}(\mathcal{S}) $.
Then, by \cite[Theorem 2.1]{CF}, we have:

\begin{enumerate}

\item[{\bf (a)}] ${\circlearrowleft_{\mathcal{S}}}$ is stable and of finite type, i.e.,  ${\circlearrowleft_{\mathcal{S}}}  =  \widetilde{\ {\circlearrowleft_{\mathcal{S}}} }$.

\item[{\bf (b)}]  The extended saturation $ \mathcal{S}^{\sharp}$ of  $\mathcal S$  coincides with $\mathcal{N}^{{\circlearrowleft_{\mathcal{S}}}}:= \{ g \in D[X] \mid  g \neq~0 \mbox{ and } \co_D(g)^{\circlearrowleft_{\mathcal{S}}} = D^{{\circlearrowleft_{\mathcal{S}}}} \} $ and ${\circlearrowleft_{\mathcal{S}}} ={\circlearrowleft_{\mathcal{S}^\sharp}}$.

\item[{\bf (c)}]  If $\mathcal{S}$ is extended saturated, then  $\Na(D, {\circlearrowleft_{\mathcal{S}}}) = D[X]_{\mathcal{S}}$.

\item[{\bf (d)}]  $\QMax^{\circlearrowleft_{\mathcal{S}}}(D) = \boldsymbol{\nabla}(\mathcal{S})$. In particular,  ${\circlearrowleft_{\mathcal{S}}}$ coincides with the spectral semistar operation associated to ${\boldsymbol{\nabla}(\mathcal{S})}$, i.e.,
$$
 E^{\circlearrowleft_{\mathcal{S}}} = \bigcap \{ED_P \mid P \in \boldsymbol{\nabla}(\mathcal{S})\}\,, \;\;\; \mbox{for all $E \in \overline{\boldsymbol{F}}(D)$}\,.
 $$

 \item[{\bf (e)}]  ${\circlearrowleft_{\mathcal{S}}}$ is a (semi)star operation on $D$ if and only if  $\mathcal{S}
\subseteq \mathcal{N}^{v_D} :=  \{ g \in D[X] \mid  g \neq~0 \mbox{ and } \co_D(g)^{v_D}  = D \} $  or, equivalently, if and only if $D = \bigcap \{ D_P \mid P \in \boldsymbol{\nabla}(\mathcal{S})\}$.

\item[{\bf (f)}] The map $\mathcal{S} \mapsto {\circlearrowleft_{\mathcal{S}}}$ establishes a 1-1 correspondence between the extended saturated multiplicative subsets of $D[X]$ (respectively, extended saturated multiplicative subsets of $D[X]$ contained in $\mathcal{N}^{v_D}$) and the set of the stable semistar (respectively, (semi)star) operations of finite type on $D$.

\item[{\bf (g)}]  Let $\mathcal{S}$ be an extended  saturated multiplicative set of $D[X]$. Then, $\Na(D, v_D) =D[X]_{\mathcal{S}}$ if and only if ${\mathcal{S}} = \mathcal{N}^{v_D}$.

\item[{\bf (h)}]  Let $R:= D^{\circlearrowleft_{\mathcal{S}}} $ and let $\iota: D \rightarrow R$ be the canonical embedding.  The overring $R$ is $t$-linked to $(D, {\circlearrowleft_{\mathcal{S}}})$ and  ${\mathcal{S}} \subseteq {\mathcal{N}}^{v_R} := \{ g \in R[X] \mid  g \neq~0 \mbox{ and } \co_R(g)^{v_R}  = R \} $ (i.e., $({\circlearrowleft_{\mathcal{S}}})_\iota$ is a (semi)star operation on $R$). Moreover $({\circlearrowleft_{\mathcal{S}}})_\iota = w_R$ if and only if  the extended saturation $
{\mathcal{S}}^{\sharp_R} := R[X] \setminus  \bigcup \{ Q[X] \mid  Q \in \Spec(R) \mbox{ and }  Q[X] \cap \mathcal S = \emptyset \}$ of the multiplicative set $\mathcal{S}$ in $R[X]$ coincides with $\mathcal{N}^{v_R}$.

\end{enumerate}
\end{remark}

 \medskip

\rm  From  Remark \ref{CF}(f),  we deduce that each semistar operation $\star$ on $D$ which is stable and of finite type is equal to
  $(\ba_{\mathcal{S}})_0 = {\circlearrowleft_{\mathcal{S}}} $
 for a unique  extended  saturated multiplicative set $\mathcal{S}$ of $D[X]$.  More precisely,

 \begin{corollary} \label{stt}
  Let $\star$ be a finite type stable semistar operation on an integral domain  $D$ with field of quotients $K$ and let $X$ be an indeterminate over $K$. Let $\mathcal{S}(\star) := \mathcal{N}^\star := \{ g \in D[X] \mid \co_D(g)^\star = D^\star \}$. Then,
 \begin{enumerate}
 \item[{\bf (1)}] $\mathcal{S}(\star)$ is an extended  saturated multiplicative set of $D[X]$ and, more precisely, $\mathcal{S}(\star) = D[X]  \setminus  \bigcup \{ Q[X] \mid  Q \in \QMax^{\stf}(D) \}$ with $\mathcal{S}(\star)^\sharp = D[X]  \setminus  \bigcup \{ P[X] \mid  P \in \Spec(D) \mbox{ and }  P[X] \cap \mathcal S(\star) = \emptyset \} = \mathcal{S}(\star) $.
 \item[{\bf (2)}] $\star = {\circlearrowleft_{\mathcal{S}(\star)}} =(\ba_{\mathcal{S}(\star)})_0$ and ${\mathcal{S}(\star)}$ is unique among the  extended  saturated multiplicative set $\mathcal{S}$ of $D[X]$  for which $ {\circlearrowleft_{\mathcal{S}}} = \star$.
 \end{enumerate}
 \end{corollary}

 \begin{proof} {\bf (1)} Clearly, $\mathcal{S}(\star) = D[X]  \setminus  \bigcup \{ Q[X] \mid  Q \in \QMax^{\stf}(D) \}$ since, for $0 \neq g \in D[X]$,  $ \co_D(g)^\star =\co_D(g)^{\stf}= D^\star$ if and only if
 $ \co_D(g)\not\subseteq Q$ for all $Q \in \QMax^{\stf}(D) $ and, for each prime ideal $Q$ of $D$, $ \co_D(g)\not\subseteq Q$ if and only if $g \not\in Q[X]$.  Moreover, clearly,
 $\bigcup \{ Q[X] \mid  Q \in \QMax^{\stf}(D) \} \subseteq \bigcup \{ P[X] \mid  P \in \Spec(D) \mbox{ and }  P[X] \cap \mathcal S(\star) = \emptyset \}$. On the other hand,
 if $P \in \Spec(D)$ and $ P[X] \cap \mathcal S(\star) = \emptyset$, this means that
 $ P[X] \subseteq \bigcup \{ Q[X] \mid  Q \in \QMax^{\stf}(D) \}$ and so $\bigcup \{ P[X] \mid  P \in \Spec(D) \mbox{ and }  P[X] \cap \mathcal S(\star) = \emptyset \} \subseteq\bigcup \{ Q[X] \mid  Q \in \QMax^{\stf}(D) \}$.

 {\bf (2)}  Note that, by assumption, $\star = \stt$ and so, for each $E \in\FF(D)$, $E^{\stt} = E\Na(D, \star) \cap K = ED[X]_{\mathcal{N}^\star} \cap K = ED[X]_{\mathcal{S}(\star)} \cap K$ \cite[Proposition 3.4]{FL}. Therefore, if $\ba_{\mathcal{S}(\star)} $ is the finite type stable semistar operation on $D[X]$ defined by the flat overring $D[X]_{\mathcal{S}(\star)} $, i.e., $A^{\ba_{\mathcal{S}(\star)}} := A_{\mathcal{S}(\star)} = AD[X]_{\mathcal{S}(\star)}$ for all $A \in \FF(D[X])$, then $\star = (\ba_{\mathcal{S}(\star)} )_0 = {\circlearrowleft_{\mathcal{S}(\star)}}$. The uniqueness follows from Remark \ref{CF}(f).
 \end{proof}

  \smallskip

Note that from  $E^{\bo}  = (E[X])^{\ba} \cap K$, by tensoring with the $D$-algebra $ D[X]$, we have $E^{\bo}[X] = (E[X])^{\ba} \cap K[X]$,
 for all $E \in \Fb(D)$. Moreover, it may happen that $E^{\bo}[X] \subsetneq (E[X])^{\ba} $
 for some $E \in \Fb(D)$.
 For instance, this happens if $E^{\bo} =K$ and if  $K[X]$ is not a $\ba$--overring of $D[X]$ (i.e., if $K[X] \subsetneq K[X]^{\ba}$).
  An explicit example is given by $\ba = v_{D[X]}$; in this case  $K^{(v_{D[X]})_0}[X] =K[X] \subsetneq (K[X])^{v_{D[X]}} =K(X)$.   Another example (even in  case of finite type stable semistar operations) is given next by  using Corollary \ref{stt}.

  \begin{example} \rm Let $P$ be a given nonzero prime   ideal of an integral domain $D$.
  Let $\Delta := \{P\}$ and set $\star := \star_\Delta$, i.e., $\star$ is the finite type stable semistar operation defined by $E^\star := ED_P$, for all $E \in \FF(D)$.  Clearly, $\QMax^{\star_\Delta}(D) = \{P\}$.  Thus, $\mathcal S(\star):= D[X] \setminus    P[X]$
  and
  $\star =\stt ={\circlearrowleft_{\mathcal S(\star)}} =   (\ba_{\mathcal S(\star)})_0 $,
   by Corollary \ref{stt} (2). (Note that $A^{\ba_{\mathcal S(\star)}} = AD[X]_{P[X]} = AD_P(X)$ for each $A \in \FF(D[X])$.)
On the other hand, for each $E \in \F(D)$, $E^\star[X] =ED_P[X] \subsetneq E D[X]_{P[X]} = E[X]D_P(X)= (E[X])^{\ba_{\mathcal S(\star)}}$
  (even if  $E^\star[X] =ED_P[X] = (E D_P(X) \cap K)[X] = E[X] D_P(X) \cap K[X] =  (E[X])^{\ba_{\mathcal S(\star)}} \cap K[X]$).
   \end{example}

  \medskip
 In order to better investigate this situation, we introduce the following definitions.
 A semistar operation $\ba$   on the polynomial domain $D[X]$ is called an {\it {extension}} (respectively, a {\it{strict extension}}) {\it{of a semistar operation}} $\star$ defined on $D$ if
$ E^{\star} = (E[X])^{\ba} \cap K$ (respectively, $ E^{\star}[X] = (E[X])^{\ba}$) for all $E \in \Fb(D)$.  Clearly, a strict extension  is an extension.  By  Lemma \ref{ast-zero},  a semistar operation $\ba$ on $D[X]$  is an extension of $\star := \bo$.

Given two semistar operations $\ba^\prime$ and $\ba^{\prime\prime}$ on the polynomial domain $D[X]$, we say that they are \it equivalent over $D$, \rm for short $\ba^{\prime}\thicksim \ba^{\prime\prime}$, (respectively, \emph{strictly equivalent  over $D$},   for short $\ba^{\prime}\thickapprox \ba^{\prime\prime}$) if   $(E[X])^{\ba^\prime} \cap K = (E[X])^{\ba^{\prime\prime}} \cap K$  (respectively, $(E[X])^{\ba^\prime} = (E[X])^{\ba^{\prime\prime}}$) for each $E \in\Fb(D)$.

Clearly, two extensions (respectively, strict extensions)  $\ba^{\prime}$ and $ \ba^{\prime\prime}$ on $D[X]$ of the same semistar operation defined on $D$ are equivalent (respectively, strictly equivalent). In particular, we have:
$$
\ba^{\prime}\thickapprox \ba^{\prime\prime} \; \Rightarrow \;   \ba^{\prime} \thicksim  \ba^{\prime\prime}\; \Leftrightarrow\; \bo^{\prime}=\bo^{\prime\prime}\,.
$$
We will see that the converse of the first implication above does not hold in general.
In order to construct some counterexamples, we need a deeper study of the problem of ``raising"  semistar operations from $D$ to $D[X]$; i.e.,
given a semistar operation $\star$ on $D$,   finding all the semistar operations $\ba$ on $D[X]$ such that $\star =  \bo$.

\smallskip

Recall that, given a family of semistar operations $\{\star_\lambda\mid \lambda \in \Lambda\}$
on an integral domain $D$, the semistar operation $\wedge  \star_\lambda$ on $D$
 is defined for all $E \in \Fb(D)$ by setting:
$$
E^{\wedge  \star_\lambda} := \bigcap  \{ E^{\star_\lambda} \mid \lambda \in \Lambda \}\,.
$$

The following statement  is a straightforward consequence of the definitions.

\begin{proposition}
\begin{enumerate}
\item[\rm\bf (1)\it ]
 Let $\star$ be a semistar operation on an integral domain $D$. Given a family of  semistar operations $\{\ba_\lambda \mid \lambda \in \Lambda\}$ on $D[X]$ that are extensions (respectively, strict extensions) of $\star$, then  $\wedge  \ba_\lambda$ is also an extension (respectively, a strict extension) of $\star$.
\item[\rm\bf (2)\it ]
 Given a family of semistar operations $\{\ba_\lambda\mid \lambda \in \Lambda\}$ of $D[X]$, suppose that $\ba_{\lambda^{\prime}}\!  \thicksim \ba_{\lambda^{\prime\prime}}$ (respectively, $\ba_{\lambda^{\prime}}\! \thickapprox \ba_{\lambda^{\prime\prime}}$) for all $\lambda^{\prime}, \lambda^{\prime\prime} \in \Lambda$,  then $\wedge  \ba_\lambda$ is equivalent  (respectively, strictly equivalent)  to $\ba_{\overline{\lambda}}$ for each $\overline{\lambda}\in \Lambda$.
\end{enumerate}
\end{proposition}

From the previous proposition, we deduce that, if a semistar operation on $D$ admits  an extension (respectively, a strict extension)    to  $D[X]$, then it admits a unique minimal extension (respectively, a unique minimal strict extension).

\medskip

At this point, it is natural to ask the following questions:

\begin{enumerate}
\item[{\bf (Q1)}]  {\sl  Given a semistar operation $\star$ defined on $D$, is it possible to find ``in a canonical way'' an extension (respectively, a strict extension) of $\star$ on $D[X]$÷?}

\item[{\bf (Q2)}]  {\sl Given an extension $\ba$ on $D[X]$ of a semistar operation $\star$ defined on $D$. Is it possible to define a strict extension   $\ba^\prime$ on $D[X]$ of $\star$ (and thus $\ba^\prime \thicksim \ba$) ?}\\
\noindent (In the statement of the previous question, we do not require  that  $\ba^\prime \thickapprox \ba$, since this condition would imply that  the extension $\ba$ on $D[X]$ was already a strict extension of $\star$.)

\end{enumerate}

\smallskip

\medskip

In the remaining part of this paper, we start the investigation of  questions {\bf (Q1)} and {\bf (Q2)}, by considering  semistar operations on $D$ defined by families of overrings. In this particular, but rather important   setting, we will provide positive answers to both questions.

\medskip

Let
$\TT := \{T_{\lambda}  \mid \lambda \in \Lambda \}$ be a nonempty set of overrings of $D$, and let $E^{{\wedge}_{\TT}}: =
\bigcap_{\lambda}ET_{\lambda}$ for each $E \in \FF(D)$. Then $\wedge_{\TT}$
is a semistar operation on $D$, and $\wedge_{\TT}$ is (semi)star if and only if
$D = \bigcap_{\lambda}T_{\lambda}$. It is easy to see that, for each $E \in \Fb(D)$ and for each $\lambda \in \Lambda$,
$$
E^{\wedge_{\TT}}T_\lambda = ET_\lambda\,,$$
(see \cite[Theorem 2]{A} for   further  details in the star operation case).   If $\TT = \{K\}$ (respectively, $\{D\}$) then obviously ${\wedge_{\{K\}}}$ (respectively, ${\wedge_{\{D\}}}$) is the trivial semistar operation $e_D \ (= \star_{\{K\}})$ (respectively, the identity (semi)star operation $d_D \ (= \star_{\{D\}})$). In case $\TT = \emptyset$, we also set ${\wedge_{\emptyset}}:= e_D$.

Note that, for each $T_{\lambda}$,
$ET_{\lambda} = \bigcap\{(ET_{\lambda})_M  \mid M \in \Max(T_\lambda)\}= \bigcap\{E(T_{\lambda})_M  \mid M \in \Max(T_\lambda)\}$; hence
$E^{\wedge_{\TT}} =
 \bigcap_\lambda \left(\bigcap\{ E(T_{\lambda})_M \mid  M  \in \Max(T_\lambda)\}\right)$.    If $\TT$ is nonempty, replacing  the family $\TT = \{T_{\lambda}\mid \lambda \in \Lambda \}$ with the family  $  \{(T_{\lambda})_M \mid  \lambda \in \Lambda, \ M \in \Max(D_\lambda)\}$, without loss of generality, whenever convenient for the context, we can assume   that each
$T_{\lambda}$ in the family of overrings $\TT$ is a quasi-local domain.

If $\TT'$ and $\TT''$ are two families of overrings of $D$, then clearly:
$$ \wedge_{\TT'} {{\bigwedge}}  \wedge_{\TT''} = \wedge_{\TT' \cup \TT''}\,.$$

Let $\TT = \{T_{\lambda} \mid \lambda \in \Lambda\}$ be a family of overrings of an integral domain $D$ with quotient field $K$.  Let $X$ be an indeterminate over $K$ and  denote by $T_{\lambda}(X)$ the Nagata ring of $T_{\lambda}$. For
each  $A \in \Fb(D[X])$, we set:
 $$ \begin{array}{rl}
 A^{\boldsymbol{\boldsymbol{(\wedge_{\TT})}}} :=&  \bigcap_{\lambda} AT_{\lambda}(X)\,,\\
 A^{\boldsymbol{{\langle\wedge_{\TT}\rangle}}} := &A^{\boldsymbol{(\wedge_{\TT})}} \cap AK[X]\,, \\
A^{\boldsymbol{[\boldsymbol{\boldsymbol{\wedge_{\TT}]}}}} :=&  \cap_{\lambda} AT_{\lambda}[X]\,.
\end{array}
$$
Clearly,
$A^{\boldsymbol{[\boldsymbol{\boldsymbol{\wedge_{\TT}]}}}} \subseteq
A^{\boldsymbol{{\langle\wedge_{\TT}\rangle}}} \subseteq
A^{\boldsymbol{\boldsymbol{(\wedge_{\TT})}}}$ for all $A \in\Fb(D[X])$, hence  ${\boldsymbol{[\boldsymbol{\boldsymbol{\wedge_{\TT}]}}}} \leq
{\boldsymbol{{\langle\wedge_{\TT}\rangle}}} \leq
{\boldsymbol{\boldsymbol{(\wedge_{\TT})}}}$.  Moreover, if $\TT$ is nonempty,
$(D[X])^{\boldsymbol{{\langle\wedge_{\TT}\rangle}}} \subseteq K[X]$, but $1/(1+X) \in  \bigcap_\lambda T_\lambda(X) = (D[X])^{\boldsymbol{\boldsymbol{(\wedge_{\TT})}}}$.  Hence,
$(D[X])^{\boldsymbol{{\langle\wedge_{\TT}\rangle}}} \subsetneq (D[X])^{\boldsymbol{\boldsymbol{(\wedge_{\TT})}}}$   and so
${\boldsymbol{{\langle\wedge_{\TT}\rangle}}} \lneq
{\boldsymbol{\boldsymbol{(\wedge_{\TT})}}}$.

\begin{proposition} \label{ext-lambda}
Let $\TT = \{(T_{\lambda}, M_\lambda)\}$ be a family of overrings of an integral domain $D$ with quotient field $K$.
Set $\boldsymbol{(\TT)} := \{T_{\lambda}(X)\mid \lambda \in \Lambda\}$, $\boldsymbol{\langle\TT\rangle} := \{T_{\lambda}(X)\mid \lambda \in \Lambda\} \cup \{K[X]\}$, and
$\boldsymbol{[\TT]} := \{T_{\lambda}[X]\mid \lambda \in \Lambda\}$. Then,
\begin{enumerate}
\item[\bf (1)] $\boldsymbol{(\wedge_{\TT})}$, $\boldsymbol{\langle\wedge_{\TT}\rangle}$, and ${\boldsymbol{[\boldsymbol{\boldsymbol{\wedge_{\TT}]}}}}$ are semistar operations on $D[X]$ and, more precisely,
$$ \boldsymbol{(\wedge_{\TT})} = \wedge_{\boldsymbol{(\TT)}},\;\;  \boldsymbol{\langle\wedge_{\TT}\rangle}= \wedge_{\boldsymbol{\langle\TT\rangle}  },\;\; \mbox{ and } \;\; \boldsymbol{[\boldsymbol{\boldsymbol{\wedge_{\TT}]}}} = \wedge_{\boldsymbol{[\TT]}} \,.$$
Moreover, if we consider the   set $\UU := \{K[X]\}$ consisting of the unique overring $K[X]$ of $D[X]$, then:
$$ \boldsymbol{\langle\wedge_{\TT}\rangle}   = \boldsymbol{(\wedge_{\TT})} \bigwedge \wedge_{\UU}  \ (=  \wedge_{\boldsymbol{(\TT)}} \bigwedge \wedge_{\UU} )\,.$$

\item[\bf (2)] The following are equivalent:
{
\begin{enumerate}
 \item[(i)] $\wedge_{\TT}$ is a (semi)star operation on $D$,
 \item[(ii)] $\boldsymbol{\langle\wedge_{\TT}\rangle}$  is a (semi)star operation on $D[X]$,
  \item[(iii)]$\boldsymbol{[\boldsymbol{\wedge_{\TT}]}}$  is a (semi)star operation on $D[X]$.
  \end{enumerate}
  }
\item[\bf (3)]   If $\TT$ is not empty and  $\TT \neq \{K\}$, then
$\boldsymbol{[\boldsymbol{\boldsymbol{\wedge_{\TT}]}}}\lneq \boldsymbol{\langle\wedge_{\TT}\rangle}$.

\item[\bf (4)] For each $E \in \Fb(D)$,
$$(E[X])^{{\bsq}} = E^{\wedge_{\TT}}[X] = (E[X])^{{\bag}} \; \mbox{ and }  \; (E[X])^{\bto} = E^{\wedge_{\TT}}(X)\,.$$
\item[\bf (5)] For each $E \in \Fb(D)$,
 $$ E^{\wedge_{\TT}} =  (E[X])^{\bsq}  \cap K = (E[X])^{\bag} \cap K =  (E[X])^{\bto} \cap K \,.$$

 \item[\bf (6)]  If $\TT$ is a finite family of overrings of $D$, then $\wedge_{\TT}$ is a semistar operation of finite type on $D$.
 \item[\bf (7)] If each overring $T \in \TT$ is a flat overring of $D$, then $\wedge_{\TT}$ is a stable semistar operation   on $D$.
\end{enumerate}
\end{proposition}
\begin{proof}
 (1) is obvious and (2) is straightforward, since it is easy to see that
 $ \bigcap \{T_{\lambda}[X]\mid \lambda \in \Lambda\} =\left(\bigcap  \{T_{\lambda}\mid \lambda \in \Lambda\}\right)[X]$ and $ \bigcap \{T_{\lambda}(X)\mid \lambda \in \Lambda\} =\left(\bigcap  \{T_{\lambda}\mid \lambda \in \Lambda\}\right)(X)$.

(3)
Without loss of generality, we can assume that $T_\lambda$ is local with nonzero  maximal ideal $M_\lambda$ and we can take the (maximal) ideal $(M_{\lambda}, 1+X)$
of $T_{\lambda}[X]$. Set $Q := (M_{\lambda}, 1+X) \cap D[X]$. Then,
$Q^{\bsq} \subsetneq Q^{\bag}$ (for
$Q^{\bsq} \cap D[X] \subseteq (M_{\lambda}, 1+X) \cap D[X] =Q
\subsetneq D[X]$ and, on the other hand, since  $M_\lambda \neq (0)$ and $Q^{\bto} = D[X]^{\bto}$, then
$Q^{\bag} \cap D[X] = Q^{\bto} \cap QK[X] \cap D[X]
= D[X]^{\bto} \cap D[X] = D[X]$).

(4) is straightforward and (5) is a trivial consequence of (4).

(6)  Clearly, for each $T \in \TT$,
 the operation $\wedge_{\{T\}}$  (defined by $A^{\wedge_{\{T\}}} := AT$,  for all $A \in \Fb(D)$, i.e., ${\wedge_{\{T\}}}= {\star_{\{T\}}}$) is a semistar operation of finite type on $D$.
  Let $\TT := \{ T_1, T_2, ...,T_n\}$ and let $z \in A^{\wedge_{\TT}}$.
  Then, $z \in F_i T_i$, for some $F_i \subseteq A$, $F_i \in \f(D)$ and for each $1 \leq i\leq n$.
    Set $F :=F_1 + F_2+ ...+ F_n$, clearly $F \in \f(D)$ and $F \subseteq A$ and $z \in FT_i$ for all $i$, thus $z \in F^{\wedge_{\TT}}$.

    (7) Let $A, B \in \Fb(D)$.  Since $T_\lambda$ is flat, $(A\cap B)T_{\lambda} = AT_{\lambda}\cap BT_{\lambda}$ \cite[Chapter 1, \S 2, N. 6]{BAC}. The conclusion is now straightforward. \end{proof}

From Proposition \ref{ext-lambda}, we deduce immediately the following:

\begin{corollary} \label{ext-sext} With the notation of Proposition \ref{ext-lambda}, assume that  $\TT$ is nonempty and $\TT \neq \{K\}$, then we have:

\begin{enumerate}

\item[\bf (1)] $\boldsymbol{[\boldsymbol{\boldsymbol{\wedge_{\TT}]}}}$, $ \boldsymbol{\langle\wedge_{\TT}\rangle}$, and  $ \boldsymbol{(\wedge_{\TT})} $ (respectively, $\boldsymbol{[\boldsymbol{\wedge_{\TT}]}}$ and $\boldsymbol{ \langle\wedge_{\TT}\rangle}$) are distinct extensions (respectively, distinct strict extensions) of $\wedge_{\TT}$.
\item[\bf (2)] $\boldsymbol{[\boldsymbol{\boldsymbol{\wedge_{\TT}]}}} \thicksim  \boldsymbol{\langle\wedge_{\TT}\rangle} \thicksim \boldsymbol{(\wedge_{\TT})} $ and, moreover, $\boldsymbol{[\boldsymbol{\wedge_{\TT}]}} \thickapprox  \boldsymbol{\langle\wedge_{\TT}\rangle}$, but neither $\boldsymbol{[\boldsymbol{\wedge_{\TT}]}}$ nor $  \boldsymbol{\langle\wedge_{\TT}\rangle}$ are strictly equivalent to $ \boldsymbol{(\wedge_{\TT})} $.

\end{enumerate}
\end{corollary}

\begin{example} \rm
 Let $\WW:= \{ W_\lambda \mid \lambda \in \Lambda\}$ be
  a family of valuation overrings of $D$ and let $\wedge_{\WW} $ be the  \texttt{ab} semistar operation on $D$ defined by the family of valuation overrings $\WW $ of $D$ (i.e., $E^{\wedge_{\WW} } := \bigcap \{ EW \mid W \in \WW \}$ for all $E \in \FF(D)$) \cite[Proposition 3.7(1)]{FL1}.
  In this case,
  \begin{enumerate}
  \item[{\bf (a)}] $\btow$ is  an \texttt{ab} semistar operation on $D[X]$ defined by the family of valuation overrings  $\WW(X):= \{ W_{\lambda}(X) \mid \lambda \in \Lambda\}$ of $D[X]$ and (by \cite[Corollary 3.8]{FL1})
for
each  $A \in \Fb(D[X])$,
$$ A^{\btow} =  \bigcap_{\lambda} AW_{\lambda}(X);$$

 \item[{\bf (b)}]  for
each  $E \in \Fb(D)$,
$$
\begin{array}{rl}
 E^{\btow_0} = & \hskip -6 pt (E[X])^{\btow} \cap K = (\bigcap_{\lambda} EW_{\lambda}(X)) \cap K \\
 = & \hskip -6 pt  \bigcap_{\lambda} (EW_{\lambda}(X) \cap K)
  = \bigcap_{\lambda} EW_{\lambda} \\
   = &\hskip -6 pt  E^{\wedge_{\WW}};
   \end{array}
   $$
    \item[{\bf (c)}]  for
each  $F \in \f(D)$,
$$
 F^{\btow_0 } =   F^{\wedge_{\WW}} =   F\Kr(D, \wedge_{\WW}) \cap K.
 $$
 In particular, by (b),  Lemma \ref{ast-zero}(2)  and \cite[Proposition 9]{FL-09}
 $$\btow_{0,f}= \btow_{f,0} =  \btow_{a,0} ={\wedge_{\WW, a}}= \wedge_{\WW, f}
 $$
 (where $\btow_{0,f} $ (respectively,
  $\btow_{f,0}$;\ $\btow_{a,0}$ ;\  ${\wedge_{\WW, a}}$;\ $ \wedge_{\WW, f}$)
 denotes the semistar operation of finite type on $D$ associated to ${\btow}_{0}$
 (respectively, the semistar operation on $D$ canonically induced by $\btow_{_{\!}f}$;\  the semistar operation on $D$ canonically induced by $\btow_{_{\!}a}$;\
the \texttt{ab} semistar operation on $D$ associated to $\wedge_{\WW}$;\  the semistar operation of finite type on $D$ associated to $\wedge_{\WW}$).
 \end{enumerate}
 \end{example}

 \begin{example}\label{d} \rm
  {\bf (1)} The identity (semi)star operation $d_D$   on an integral domain $D$, is defined by the family of a single overring $\DD := \{D\}$ of $D$, i.e., $d_D = \wedge_{\DD}$. Set
$$
\boldsymbol{[}d_D\boldsymbol{]}:=\bsqd\,,\;\;
\boldsymbol{\langle} d_D \boldsymbol{\rangle} := \bagd \,,\;\;  \boldsymbol{(} d_D \boldsymbol{)} := \btod\,. $$
Clearly,   if $D$ is not a field, $d_{D[X]} =\boldsymbol{[}d_D\boldsymbol{]} \lneq
\boldsymbol{\langle} d_D \boldsymbol{\rangle} \lneq \boldsymbol{(} d_D \boldsymbol{)}$ and  $\boldsymbol{\langle} d_D \boldsymbol{\rangle} $ (and  $\boldsymbol{[}d_D\boldsymbol{]}$) is a (semi)star operation on $D[X]$, but in general $\boldsymbol{(} d_D \boldsymbol{)}$ is not  a (semi)star operation on $D[X]$.
Moreover,  $\boldsymbol{\langle} d_D \boldsymbol{\rangle} $, $\boldsymbol{(} d_D \boldsymbol{)}$ (and  $\boldsymbol{[}d_D\boldsymbol{]}$) are stable semistar operations of finite type, since  $\boldsymbol{\langle} d_D \boldsymbol{\rangle} $ is defined by the two flat overrings $D(X)$ and $K[X]$ of $D[X]$ and  $\boldsymbol{(} d_D \boldsymbol{)}$ is defined by a unique flat overring $D(X)$ of $D[X]$ (Proposition \ref{ext-lambda}((6) and (7))).

 {\bf (2)} As observed above, the trivial semistar operation $e_D$ on an integral domain $D$, with quotient field $K$, is defined by the family of a single overring $\KK:=\{K\}$ of $D$, i.e., $e_D = \wedge_{\KK}$. Set
$$
\boldsymbol{[}e_D\boldsymbol{]}:=\bsqk\,,\;\;
\boldsymbol{\langle} e_D \boldsymbol{\rangle} := \bagk \,,\;\;  \boldsymbol{(} e_D \boldsymbol{)} := \btok\,. $$

Clearly, $\boldsymbol{[}e_D\boldsymbol{]} =
\boldsymbol{\langle} e_D \boldsymbol{\rangle} \lneq \boldsymbol{(} e_D \boldsymbol{)} =e_{D[X]}$, where   $\boldsymbol{[}e_D\boldsymbol{]} \ (=\boldsymbol{\langle} e_D \boldsymbol{\rangle})$  is the stable semistar operation of finite type on $D[X]$
 defined by the  flat overring $K[X]$, i.e., $\boldsymbol{[}e_D\boldsymbol{]} =\boldsymbol{\langle} e_D \boldsymbol{\rangle} = {\boldsymbol \wedge}_{\{K[X]\}} \ (= {\boldsymbol \star}_{\{K[X]\}})$.
\end{example}
 \medskip

We study now the important case in which the family of valuation overrings $\WW$ of $D$ coincides with the family of all valuation overrings of $D$.

\begin{proposition}\label{b} Let $\VV $ be the family of all valuation overrings of an integral domain $D$ with quotient field $K$.
Note that $\wedge_{\VV}$ coincides with $b_D$ (the $b$--operation on $D$; see Section 1).  Set
$$
\boldsymbol{[}b_D\boldsymbol{]}:=\bsqv\,,\;\;
\boldsymbol{\langle} b_D \boldsymbol{\rangle} := \bagv \,,\;\;  \boldsymbol{(} b_D \boldsymbol{)} := \btov\,. $$

\begin{enumerate}
\item[{\bf (1)}] $\boldsymbol{[}b_D\boldsymbol{]}$, $
\boldsymbol{\langle} b_D \boldsymbol{\rangle} $, and $\boldsymbol{(} b_D \boldsymbol{)}$ are semistar operations on $D[X]$ with $\boldsymbol{[}b_D\boldsymbol{]} \leq
\boldsymbol{\langle} b_D \boldsymbol{\rangle} \leq\boldsymbol{(} b_D \boldsymbol{)}$.  If $D$ is integrally closed then $\boldsymbol{[}b_D\boldsymbol{]}$ and $
\boldsymbol{\langle} b_D \boldsymbol{\rangle} $ are (semi)star operations on $D[X]$.
In general, $\boldsymbol{(} b_D \boldsymbol{)}$ is not a (semi)star operation on $D[X]$ even if $D$ is integrally closed (or, equivalently, even if $b_D$ is a (semi)star operation on $D$).

\item[{\bf (2)}]  If $D\neq K$, i.e. if $D$ has at least one nontrivial valuation overring,  then
$\boldsymbol{[}b_D\boldsymbol{]}$, $
\boldsymbol{\langle} b_D \boldsymbol{\rangle} $, and $\boldsymbol{(} b_D \boldsymbol{)}$ (respectively, $\boldsymbol{[}b_D\boldsymbol{]}$ and  $
\boldsymbol{\langle} b_D \boldsymbol{\rangle} $) are distinct extensions (respectively, distinct strict extensions) of $b_D$.

\item[{\bf (3)}]
$\boldsymbol{\langle} b_D \boldsymbol{\rangle}$ and  $\boldsymbol{(} b_D \boldsymbol{)}$ are \texttt{ab} semistar operations such that
  $$
  \boldsymbol{[}b_D\boldsymbol{]} \leq \boldsymbol{[}b_D\boldsymbol{]}_a = b_{D[X]} \leq \boldsymbol{\langle} b_D \boldsymbol{\rangle} \leq \boldsymbol{(} b_D \boldsymbol{)}
  $$
and, in general, $\boldsymbol{[}b_D\boldsymbol{]}$
is not an \texttt{eab} semistar operation.
\end{enumerate}
\end{proposition}
\begin{proof} (1)  and (2) follow from Proposition \ref{ext-lambda}((1), (2) and (3)) and Corollary \ref{ext-sext}.

(3) Clearly, $ \boldsymbol{(} b_D \boldsymbol{)}$ is an \texttt{ab}
 semistar operation on $D[X]$, since it is defined by the family of valuation overrings $\VV(X) :=\{V(X) \mid V \in \VV\}$.

 Moreover, if \  $\UU := \{K[X]\}$ \  and \ $\UU' := \{K[X]_M \mid  M\in \Max(K[X])\}$, then clearly $\wedge_{\UU} = \wedge_{\UU'}$. Therefore,  $\boldsymbol{\langle} b_D \boldsymbol{\rangle} = \boldsymbol{\langle\wedge_{\VV}\rangle}   = \boldsymbol{(\wedge_{\VV})} \bigwedge \wedge_{\UU}  =  \wedge_{\VV(X)} \bigwedge \wedge_{\UU'}$ and, hence $\boldsymbol{\langle} b_D \boldsymbol{\rangle}$ is also an \texttt{ab}
 semistar operation on $D[X]$.

  Note that $ b_{D[X]} \leq \boldsymbol{\langle} b_D \boldsymbol{\rangle}$ because \  $\VV(X) \bigcup \UU' $ \ is a subset of the family of all valuation overrings of $D[X]$.

 Let $W$ be a valuation overring of $D[X]$ with maximal ideal $N$. Two cases are possible. If $W \cap K = K$, then $K[X] \subseteq W$,
and hence $W = K[X]_M$ for some maximal ideal $M$ of $K[X]$. Next, if $W \cap K \subsetneq K$,
then $V:=W \cap K$ is a valuation overring of $D$ with nonzero maximal ideal $\boldsymbol{n}$,
and so $W$ is a valuation overring of $V[X]$.  If $W \neq   V(X) =V[X]_{\boldsymbol{n}[X]}$ the maximal ideal $N$ of $W$ must contract on a maximal ideal of $V[X]$ upper to the maximal ideal $\boldsymbol{n}$ of $V$, i.e., $N \cap V[X] \supsetneq \boldsymbol{n}[X]$.  Therefore, we have necessarily that $V[X] \subset W \subseteq V(X)$.
From the previous observations, we easily deduce that $\boldsymbol{[}b_D\boldsymbol{]}_a = b_{D[X]}$.


We next construct an integral domain $D$ such that
$\boldsymbol{[}b_D\boldsymbol{]}$ is not an \texttt{eab} semistar operation.
 Let $D:= \mathbb R + T \mathbb C [\![T]\!]$, i.e., $D$ is a pseudo-valuation domain with canonically associated valuation overring $V:=\mathbb C [\![T]\!]$ and quotient field $K := \mathbb C(\!(T)\!)$.  Since $\mathbb R \subset \mathbb C$ is a finite field extension the valuation overrings of $D$ are just $V$ and $K$, thus it is straightforward to see that  that $d_D \lneq b_D = \wedge_{\{V\}}$  and
$$
d_{D[X]} \lneq \boldsymbol{[}b_D\boldsymbol{]} =\wedge_{\{V[X],\ K[X]\}} \leq   b_{D[X]} \leq \boldsymbol{\langle} b_D \boldsymbol{\rangle} = \wedge_{\{V(X),\ K[X]\}} \lneq  \boldsymbol{(} b_D \boldsymbol{)} = \wedge_{\{V(X)\}}$$
 (Proposition \ref{ext-lambda}(3)).
 Moreover, $[b_D]$ is not an \texttt{eab} semistar operation on $D[X]$, because
 if $\boldsymbol{[}b_D\boldsymbol{]} \ (= \wedge_{\{V[X],\ K[X]\}}= \star_{\{V[X]\}})$ was an \texttt{eab} semistar operation on $D[X]$, since it is of finite type, then 
 $b_{D[X]} = (d_{D[X]})_a \leq  
 \boldsymbol{[}b_D\boldsymbol{]}  \leq 
  b_{D[X]}$, i.e. $ \star_{\{V[X]\}} = \boldsymbol{[}b_D\boldsymbol{]} = b_{D[X]}$, which is a contradiction since $V[X]$ is not a Pr\"ufer domain.
\end{proof}

\medskip

In the next result, we provide another application of Proposition \ref{ext-lambda}.

\begin{proposition} \label{M}  Let $\star$ be a semistar operation of an integral domain $D$ with quotient field $K$ and let $X$ be an indeterminate over $K$. Set $\MM := \MM(\star) := \{ D_Q \mid Q \in \QMax^{\stf}(D)\}$.   It is well known that, in this case,  $\wedge_{\!\MM}$ coincides with $ \stt$, the stable semistar operation of finite type associated to $\star$. Set
$$    \boldsymbol{(}\stt\boldsymbol{)} := \btom \,,\;
 \boldsymbol{\langle}\stt\boldsymbol{\rangle}  :=  \bagm\,,\;
 \mbox{and } \;   \boldsymbol{[}\stt\boldsymbol{]} := \bsqm\,.$$
 \begin{enumerate}
\item[\bf (1)]
 For each $A \in \Fb(D[X])$,
$$
 \begin{array}{rl}
A^{\boldsymbol{[}\stt\boldsymbol{]}} = & \hskip -5pt \bigcap \{AD_Q[X] \mid Q \in \QMax^{\stf}(D) \}\,,  \\
A^{\boldsymbol{\langle}\stt\boldsymbol{\rangle}} = & \hskip -5pt  A\Na(D, \star) \cap AK[X]\,, \mbox{ and } \quad  A^{\boldsymbol{(}\stt\boldsymbol{)}} = A\Na(D, \star)\,.
\end{array}
$$
\item[\bf (2)]
 ${\boldsymbol{[}\stt\boldsymbol{]}}$,
 ${\boldsymbol{\langle}\stt\boldsymbol{\rangle}}$, and
 ${\boldsymbol{(}\stt\boldsymbol{)}}$
 are stable semistar operations of $D[X]$; moreover,
  ${\boldsymbol{\langle}\stt\boldsymbol{\rangle}}$ and
 ${\boldsymbol{(}\stt\boldsymbol{)}}$
 are also of finite type. Therefore,
$$
{\boldsymbol{\langle}\stt\boldsymbol{\rangle}} =
{\boldsymbol{\langle}\stt\boldsymbol{\rangle}}_{\!_f} =
{\widetilde{\ {\boldsymbol{\langle}\stt\boldsymbol{\rangle}}\ }},\;
{\boldsymbol{(}\stt\boldsymbol{)}} = {\boldsymbol{(}\stt\boldsymbol{)}}_{\!_f} = \widetilde{\ {\boldsymbol{(}\stt\boldsymbol{)}}\ }\,.$$
\item[\bf (3)]
 For each $E \in \Fb(D)$,
$$ (ED[X])^{\boldsymbol{[}\stt\boldsymbol{]}} = E^{ \stt}[X] =
(ED[X])^{\boldsymbol{\langle}\stt\boldsymbol{\rangle}} \quad
\mbox{ and } \quad  (ED[X])^{\boldsymbol{(}\stt\boldsymbol{)}} = E\Na(D, \star)\,;$$
and so
 $$ E^{ \stt} =  (ED[X])^{\boldsymbol{[}\stt\boldsymbol{]}} \cap K =
 (ED[X])^{\boldsymbol{\langle}\stt\boldsymbol{\rangle}} \cap K =
 (ED[X])^{\boldsymbol{(}\stt\boldsymbol{)}} \cap K \,.$$

\end{enumerate}
\end{proposition}
\begin{proof}

 {(1)} and  {(2)}. Note that  the semistar operation  ${\boldsymbol{[}\stt\boldsymbol{]}}$ (respectively,\
 ${\boldsymbol{\langle}\stt\boldsymbol{\rangle}}$;\
 ${\boldsymbol{(}\stt\boldsymbol{)}}$) on $D[X]$ is defined by the family of flat overrings $\{D_Q[X] \mid Q \in \QMax^{\stf}(D) \}$ (respectively, \ $\{D_Q(X) \mid Q \in \QMax^{\stf}(D) \} \cup \{K[X]\}$\ \!; \ $\{D_Q(X) \mid Q \in \QMax^{\stf}(D) \}$÷) of $D[X]$ and so is a stable semistar operation of $D[X]$.  The first equality in (1) is a transcription of the definition.
The last two equalities in (1)  are consequence of the fact that $\bigcap \{AD_Q(X)\mid Q \in \QMax^{\stf}(D) \} = A\Na(D, \star)$ (Proposition \ref{nagata}(2) or \cite[Proposition 3.1(3)]{FL}). Finally, it is easy to see that ${\boldsymbol{(}\stt\boldsymbol{)}}$ (respectively,\ ${\boldsymbol{\langle}\stt\boldsymbol{\rangle}}$
 ) is a (stable) semistar operation of finite type, since it is defined by a unique flat overring of $D[X]$, i.e.,  $\Na(D, \star)$ (respectively, by two
 flat overrings of $D[X]$, i.e., $\Na(D, \star)$ and $K[X]$) (Proposition \ref{ext-lambda}((6) and (7)).

  {(3)}  is an application of Proposition \ref{ext-lambda}((4) and (5)).
\end{proof}

\medskip

It is natural to ask if, eventually,  ${\boldsymbol{[}\stt\boldsymbol{]}}  \neq {\boldsymbol{[}\stt\boldsymbol{]}}_{\!{_f}}$. The answer to this question is negative (i.e., ${\boldsymbol{[}\stt\boldsymbol{]}} $ is also a semistar operation of finite type (and stable)).
In order to show this fact, we deepen the study of the semistar operation ${\boldsymbol{[}\stt\boldsymbol{]}}$   defined on $D[X]$.

\begin{proposition} \label{[star]}
Let  $D$, $X$, $\star$  and ${\boldsymbol{[}\stt\boldsymbol{]}}$ be as in Proposition \ref{M}.
 Then, for all $A \in \FF(D[X])$,
$$
\begin{array}{rl}
A^{{\boldsymbol{[}\stt\boldsymbol{]}}} &=
\bigcup\{ (A: F) \mid F \in \f(D) \mbox{ and }  F^{\stt}= D^{\stt}\}\\
&= \bigcup\{ (A: H) \mid H \in \f(D),\ H\subseteq D, \mbox{ and }  H^{\stt}= D^{\stt}\} \,.
\end{array}
$$
In particular, ${\boldsymbol{[}\stt\boldsymbol{]}}$ is a (stable) semistar operation of finite type on $D[X]$ and so
 ${\boldsymbol{[}\stt\boldsymbol{]}}=
 {\boldsymbol{[}\stt\boldsymbol{]}}_{\!{_f}} =
 \widetilde{\ {\boldsymbol{[}\stt\boldsymbol{]}}\ }$.
\end{proposition}

\begin{proof} Set $A_1 := \bigcup\{ (A: F) \mid F \in \f(D) \mbox{ and } F^{\stt}= D^{\stt}\} $ and $A_2 := \bigcup\{ (A: H) \mid H \in \f(D),\ H\subseteq D, \mbox{ and }  H^{\stt}= D^{\stt}\}$. Clearly, $A_2 \subseteq A_1$.

We start by showing that $A_1 \subseteq A^{{\boldsymbol{[}\stt\boldsymbol{]}}}$.  Note that $F^{\stt} = D^{\stt}$ if and only if $ FD_Q = D_Q$ for all $Q \in \QMax^{\stf}(D) =\QMax^{\stt}(D)$.
If $z \in A_1$, then $zF \subseteq A$ for some $ F \in \f(D)$ and $F^{\stt} = D^{\stt}$, hence $ z \in z D_Q[X] = zFD_Q[X] \subseteq AD_Q[X]$ for all $Q \in \QMax^{\stf}(D)$ and so $ z \in A^{{\boldsymbol{[}\stt\boldsymbol{]}}} $.

Now, we show that $A^{{\boldsymbol{[}\stt\boldsymbol{]}}}\subseteq A_2$. Let $z \in A^{{\boldsymbol{[}\stt\boldsymbol{]}}} = \bigcap \{ AD_Q[X] \mid Q \in \QMax^{\stf}(D)\} $ and set $I := \{d \in D \mid dz \in A \}$. It is easy to see that $I$ is an ideal of $D$ (depending on $A$ and $z$) and moreover $I \not\subseteq Q$, i.e., $ ID_Q = D_Q$ for all  $Q \in \QMax^{\stf}(D)$. Since $I^{\stt} = D^{\stt}$, we can find $H\in \f(D)$
with $H \subseteq I$ and $H^{\stt} = D^{\stt}$.  Therefore $zH \subseteq zI \subseteq A$ and so $z \in A_2$.

The last statement follows easily from Proposition \ref{M}(2) and from  the fact that, if $z \in A^{{\boldsymbol{[}\stt\boldsymbol{]}}}=\bigcup\{ (A: F) \mid F \in \f(D) \mbox{ and }  F^{\stt}= D^{\stt}\}$,
 then $zF_0 =: G_0 \subseteq A$
  for some $F_0 \in \f(D)$, with $F_0^{\stt}= D^{\stt}$, and so $G_0 \in  \f(D)$,
  and $z \in (G_0 : F_0) \subseteq G_0^{\boldsymbol{[}\stt\boldsymbol{]}} \subseteq
  A^{{\boldsymbol{[}\stt\boldsymbol{]}}_{\!{_f}}} $.
\end{proof}

From Corollaries \ref{stt} and  \ref{ext-sext} and from Propositions  \ref{M} and \ref{[star]}, we easily obtain the following:

\begin{corollary} Let $\star$, $\stt$,  $D$, $K$,  $X$, ${\boldsymbol{[}\stt\boldsymbol{]}}$,
 ${\boldsymbol{\langle}\stt\boldsymbol{\rangle}}$, and
 ${\boldsymbol{(}\stt\boldsymbol{)}}$ be as in Proposition \ref{M}.  Let ${\mathcal S(\star)}$ and $\ba_{\mathcal S(\star)}$ be as in Corollary \ref{stt}.  Assume that   $D_Q \subsetneq K$  for some $Q \in \QMax^{\stf}(D)$.
 \begin{enumerate}

\item[\bf (1)] ${\boldsymbol{[}\stt\boldsymbol{]}}$,
 ${\boldsymbol{\langle}\stt\boldsymbol{\rangle}}$, and
 ${\boldsymbol{(}\stt\boldsymbol{)}}$ (respectively, ${\boldsymbol{[}\stt\boldsymbol{]}}$,
 and ${\boldsymbol{\langle}\stt\boldsymbol{\rangle}}$) are distinct  finite type stable semistar extensions to $D[X]$ (respectively, distinct  finite type stable semistar strict extensions) of \ $\stt$.  Moreover, ${\boldsymbol{(}\stt\boldsymbol{)}}$ is the unique finite type stable semistar extension to $D[X]$ of \ $\stt$ defined by an extended saturated multiplicative set of $D[X]$, i.e.,  ${\boldsymbol{(}\stt\boldsymbol{)}} = \ba_{\mathcal S(\star)}$, where
 ${\mathcal S(\star)} = \mathcal{N}^\star := \{ g \in D[X] \mid \co_D(g)^\star = D^\star \}$.
\item[\bf (2)] ${\boldsymbol{[}\stt\boldsymbol{]}} \thicksim {\boldsymbol{\langle}\stt\boldsymbol{\rangle}} \thicksim {\boldsymbol{(}\stt\boldsymbol{)}}$ and, moreover, ${\boldsymbol{[}\stt\boldsymbol{]}} \thickapprox {\boldsymbol{\langle}\stt\boldsymbol{\rangle}}$, but neither ${\boldsymbol{[}\stt\boldsymbol{]}}$ nor ${\boldsymbol{\langle}\stt\boldsymbol{\rangle}}$ are strictly equivalent to $ {\boldsymbol{(}\stt\boldsymbol{)}} $. \hfill $\Box$
\end{enumerate}
\end{corollary}

\begin{remark} \label{Pic} \rm G. Picozza \cite{Pi} has studied a  different approach for the extension to the polynomial ring  $D[X]$ of a semistar operation defined on  an integral  domain $D$.

First, he proves the following \cite[Propositions 3.1 and 3.2]{Pi}.
\begin{enumerate}
\item[{\bf (a)}] \it Let $\cal F$ be a localizing system of ideals of $D$ and set:
$$ \cal F [X] : = \{ J  \mbox{ ideal of $D[X] \mid  J \supseteq I[X]$ for some ideal $I$ of $\cal F$}\}\,.
$$
 { 
 \begin{enumerate}
 \item[\bf{(a.1)}]  \it $\cal F [X]$ is a localizing system on $D[X]$.
 \item[\bf{(a.2)}] \it $\cal F [X] =\{J  \mbox{ ideal of $D[X] \mid J\cap D \in \cal F$}\}$.
 \item[\bf{(a.3)}] \it If $\cal F$ is a localizing system of finite type of $D$, then $\cal F[X]$ is a localizing system of finite type on $D[X]$
\end{enumerate}
} 
\end{enumerate}

 Then, he uses some of the results by Fontana and Huckaba recalled in Proposition \ref{prop:loc1}. More precisely, if $\cal F$ is a localizing system on $D$, Picozza considers the semistar operation $\star_{\cal F[X]}$ on $D[X]$ canonically associated to the localizing system $\cal F[X]$ on $D[X]$ introduced in (a). In particular, if $\cal F $ is the localizing system  associated to a given semistar operation $\star$ defined on $D$, i.e., $\cal F  = \cal F^\star$, he considers the stable semistar operation  on $D[X]$ associated to the localizing system ${\cal F^\star [X]}$. Set $\star[X] := \star_{\cal F^\star [X]}$ (Picozza denotes by $\star^\prime$ this semistar operation on $D[X]$ \cite[Theorem 3.3]{Pi}).

 We are now in a position to compare the semistar operations on the polynomial rings studied by Picozza and the semistar operation ${\boldsymbol{[}\stt\boldsymbol{]}}$ introduced in Proposition \ref{M}. The following result improves \cite[Proposition 3.4]{Pi}.
 \begin{enumerate}
 \item[\bf{(b)}]  \it Using the notation introduced above, then
 $$
 \begin{array}{rl}
 {\boldsymbol{[}\stt\boldsymbol{]}} =& \hskip -6pt  \stt[X] = {\widetilde{\ \star[X]}\ } \\
 = & \hskip -6pt  \stf[X] = \star[X]_{\!{_f}}\,.
 \end{array}
 $$
\end{enumerate}\rm

It is clear that $\stt[X] = \stf[X]$, since:
$$
\begin{array}{rl}
\cal F^{\stf} = & \hskip -5 pt \{ I \mbox{ ideal of $D  \mid I^{\stf} =D^{\stf}$}\}\\
= & \hskip -5 pt \{ I \mbox{ ideal of $D  \mid I\not\subseteq Q$ for all $Q \in \QMax^{\stf}(D)$}\}  \\
= & \hskip -5 pt \{I \mbox{ ideal of $D  \mid I \not\subseteq Q$ for all $Q \in \QMax^{\stt}(D)$} \}  \\
= & \hskip -5 pt \{ I \mbox{ ideal of $D  \mid I^{\stt} =D^{\stt}$}\} \\
= & \hskip -5 pt \cal F^{\stt}\,.
\end{array}
$$
By Proposition \ref{[star]}, we know that $ A^{\boldsymbol{[}\stt\boldsymbol{]}}  =
\bigcup\{ (A: H) \mid H\in \f(D),\ H \subseteq D \mbox{ and }  H^{\stt}= D^{\stt}\}$ for all $A \in \FF(D[X])$. On the other hand (by definition of a semistar operation associated to a localizing system (Proposition \ref{prop:loc1}(3)) and by (a.2)), we have  $ A^{{\stt}[X]} =
\bigcup\{ (A: J) \mid J \mbox{ ideal of $D[X] $ such that } $ $(J\cap D)^{\stt}= D^{\stt}\}$ for all $A \in \FF(D[X])$.  Therefore, if $z \in A^{\boldsymbol{[}\stt\boldsymbol{]}} $, then $z \in (A:H) = (A: H[X])$ for some finitely generated ideal $H$ of $D$ such that  $(H[X] \cap D)^{\stt} = D^{\stt}$, thus $z \in
A^{{\stt}[X]}$.

Conversely, let $z \in A^{{\stt}[X]} $. Therefore, $z \in  (A: J)$ for some ideal $J$ of $D[X]$ such that $(J\cap D)^{\stt}= D^{\stt}$. In this situation, we can find a finitely generated ideal $H$ in $D$ such that $H \subseteq J\cap D$ and $H^{\stt} = D^{\stt}$. Since $H \subseteq J$ then $(A:J) \subseteq (A:H)$,  thus $z \in (A:H)$ and so $z \in
A^{\boldsymbol{[}\stt\boldsymbol{]}} $.

From the previous results, we deduce:
\begin{enumerate}
 \item[\bf{(c)}]  \it $\boldsymbol{[}d_D\boldsymbol{]} =d_D[X] = d_{D[X]}$.
  \item[\bf{(d)}]  \it  $\boldsymbol{[}w_D\boldsymbol{]} =w_D[X]  =t_D[X] \leq w_{D[X]} \leq t_{D[X]}$.

   \end{enumerate} \rm

The statement (c) is a straightforward consequence of (b), since $d_D = {\widetilde{\ d_D \ }} = d_{D, f}$ and the localizing system $\cal F^{d_D}[X] = \{ D[X] \} =\cal F^{d_{D[X]}}$.

The equalities in (d) are obtained from   (b) (and Proposition \ref{M}) by taking $\star = v_D$ (and so, $\stt = w_D$ and $\stf = t_D$).
Moreover, it is always true that the $w$--operation is   smaller than  or equal to the $t$--operation.
Finally, note that $w_D[X]$ (respectively, $w_{D[X]}$) is the stable (semi)star operation of finite type on $D[X]$ canonically associated to the localizing system $\cal F^{w_D}[X] =\{ J \mbox{ ideal of } D[X] \mid  J \supseteq I[X] \mbox{ with } I^{w_D} = D \}$ (respectively,
 $\cal F^{w_{D[X]}} =\{ J \mbox{ ideal of }  D[X] \mid  J ^{w_{D[X]}} = D[X] \}$)  \cite[Theorem 2.10 (B)]{FH}.  Since $I^{w_D} = D$ implies that $I^{w_D}[X] = D[X]$ and so also $I[X]^{w_{D[X]}} =D[X]$ (see  the proof of Lemma \ref{ast-zero}(7)), then
$\cal F^{w_D}[X] \subseteq \cal F^{w_{D[X]}}$. From this, we conclude that ${w_D}[X]  \leq w_{D[X]}$.

Note that, in (d),  it may happen that  $\boldsymbol{[}w_D\boldsymbol{]} =w_D[X] \lneq w_{D[X]}$ (e.g., by \cite[Remark 2]{Pi}, let  $Q$ be a prime ideal of $D[X]$ not extended from $D$ and such that $Q\cap D$ is a  $t_D$--maximal ($= w_D$--maximal) ideal of $D$;
 since $(Q \cap D)[X] \subsetneq Q$ and $(Q \cap D)[X]$ is a
 $t_{D[X]}$--maximal ($= w_{D[X]}$--maximal) ideal of $D[X]$
 (\cite[Proposition 1.1]{HZ} and \cite[Corollary 3.5 (2)]{FL0}), then $Q$ is not a $w_{D[X]}$--(maximal) ideal, but clearly $Q$ is a $\boldsymbol{[}w_D\boldsymbol{]}$--ideal).
\end{remark}
\bigskip

 \begin{remark}\label{sahandi} \rm {\bf (1)} Using the techniques introduced in \cite{CF} and recalled in Remark \ref{CF}, Sahandi \cite[Theorem 2.1, Proposition 2.2 and Remark 2.3]{S} has recently given another
 description of the stable semistar operation ${\boldsymbol{[}\stt\boldsymbol{]}}$.

Let $D{_{\!{_1}}}:=D[X]$, $K{_{\!{_1}}}:=K(X)$, $Y$ an indeterminate over $K{_{\!{_1}}}$, and  consider the following subset of $\Spec(D{_{\!{_1}}})$:
 $$
 \Delta{_{\!{_1}}} := \Delta{_{\!{_{1, \star}}}} := \{Q{_{\!{_1}}}  \in \Spec(D{_{\!{_1}}}) \mid
 \mbox{either }   Q{_{\!{_1}}} \cap D = (0)  \mbox{ or } (Q{_{\!{_1}}} \cap D)^{\stf} \subsetneq D^{\stf} \}\,.
 $$
 Set
 $$
 \mathcal S^{\sharp}_{\!{_1}} := \mathcal S^{\sharp}_{\!{_{1, \star}}}:= \mathcal S^{\sharp}(\Delta{_{\!{_{1, \star}}}}) :=
D{_{\!{_1}}}[Y] \setminus \left(\bigcup \{ Q{_{\!{_1}}}[Y] \mid Q{_{\!{_1}}} \in \Delta{_{\!{_{1, \star}}}} \}\right)\,.
$$
Clearly, $\mathcal S^{\sharp}_{\!{_1}}$ is an extended saturated multiplicative system of
$D{_{\!{_1}}}[Y]$ and so we can consider  the stable semistar operation of finite type on
 $D{_{\!{_1}}}$,
 $\circlearrowleft_{\mathcal S^{\sharp}_{\!{_1}}}$, defined by setting for each $A \in \FF(D_{\!{_1}}) = \FF(D[X])$:
 $$
 A^{\circlearrowleft_{\mathcal S^{\sharp}_{\!{_1}}}} :=
 AD_{\!{_1}}[Y]_{\mathcal S^{\sharp}_{\!{_1}}} \cap K_{\!{_1}}\,.
 $$

 We have already observed in Remark \ref{Pic}(b) that ${\boldsymbol{[}\stt\boldsymbol{]}}$ is the unique stable semistar operation on $D_{\!{_1}} \ (= D[X])$ determined by the localizing system of finite type $ \mathcal F_{\!{_{1}}} :=  \mathcal F_{\!{_{1, \star}}} := \mathcal F^{\stf}[X] =\mathcal F^{\stt}[X]$.
 Moreover, the map  $ \boldsymbol{\sigma}: \mathcal F_{\!{_1}} \mapsto \mathcal S(\mathcal F_{\!{_1}})
 := D_{\!{_1}}[Y] \setminus  \bigcup \{ P{_{\!{_1}}}[Y] \mid P{_{\!{_1}}} \not\in  \mathcal F_{\!{_1}}\}$ establishes a bijection between the set of  the localizing systems of finite type on $D_{\!{_1}}$ and the extended saturated multiplicative systems of
$D{_{\!{_1}}}[Y]$ and, under this map, the corresponding associated stable semistar operations of finite type on $D{_{\!{_1}}}$ coincide \cite[Corollary 2.2]{CF}.  Since it is straightforward that $\Delta{_{\!{_{1, \star}}}}$ coincides with the set $ \{ P{_{\!{_1}}} \in \Spec(D_{\!{_1}}) \mid P{_{\!{_1}}} \not\in  \mathcal F_{\!{_{1, \star}}}\}$, then clearly  $\mathcal S^{\sharp}_{\!{1, \star}}$ corresponds canonically to $\mathcal F_{\!{_{1, \star}}}$ under $ \boldsymbol{\sigma}$. Therefore ${\boldsymbol{[}\stt\boldsymbol{]}}$ coincides with $\circlearrowleft_{\mathcal S^{\sharp}_{\!{_1}}}$.
Moreover, we also have $\Na(D_{\!{_1}}, {\boldsymbol{[}\stt\boldsymbol{]}}) = D_{\!{_1}}[Y]_{\mathcal S^{\sharp}_{\!{_1}}}$
 or, equivalently, ${\mathcal S^{\sharp}_{\!{_1}}} =
 \{0 \neq g_{\!{_1}}  \in D_{\!{_1}}[Y] \mid \co_{D_{\!{_1}}}(g_{\!{_1}})^{{\boldsymbol{[}\stt\boldsymbol{]}}} = D_{\!{_1}}^{{\boldsymbol{[}\stt\boldsymbol{]}}} \}=:  \cal N_{\!{_1}}^{\boldsymbol{[}\stt\boldsymbol{]}}$ \cite[Theorem 2.1((c) and (d))]{CF}.

 Let $ {\nabla}_{\!{_{1, \star}}} $ be the set of the maximal elements of ${\Delta}_{\!{_1, \star}}$.

  It is easy to see that
 $$
 \begin{array}{rl}
{\nabla}_{\!{_{1, \star}}} = \{ Q_{\!{_1}} \in \Spec(D{_{\!{_1}}}) \mid &
  \mbox{either } \;\; Q{_{\!{_1}}} \cap D = (0) \mbox{ and } \co_D(Q{_{\!{_1}}})^{\stf} = D^{\stf}  \\
  & \mbox{or } \;\; Q{_{\!{_1}}} \cap D \in \QMax^{\stf}(D) \}\,,
  \end{array}
$$
since a prime ideal $ Q{_{\!{_1}}} \in \Spec(D{_{\!{_1}}})$ such that $Q{_{\!{_1}}} \cap D = (0)$ is not contained in any ideal of the type $Q[X]$ with $Q \in \QMax^{\stf}(D)$ if and only if  $ \co_D(Q{_{\!{_1}}})^{\stf} = D^{\stf}$.

For the sake of simplicity, set $\boldsymbol{\star}_{\!{_1}} := {\circlearrowleft_{\mathcal S^{\sharp}_{\!{_1}}}} \ (=  \boldsymbol{[}\stt\boldsymbol{]})$.
Then, by the previous remarks, we can conclude that $ \QMax^{\boldsymbol{\star}_{\!{_1}}}(D_{\!{_1}}) ={\nabla}_{\!{_{1, \star}}} $ \cite[Theorem 2.1(e)]{CF}.

Putting together the previous information with Proposition \ref{M}(1), for all $A \in \FF(D[X])$, we have:
$$
\begin{array}{rl}
\bigcap \{AD_Q[X] \mid Q \in \QMax^{\stf}(D) \} = A^{\boldsymbol{[}\stt\boldsymbol{]}}
= & \hskip -5pt
\bigcap \{AD[X]_{Q{_{\!{_1}}}} \mid Q{_{\!{_1}}} \in {\nabla}_{\!{_{1, \star}}}  \} \\
= & \hskip -5pt  AD[X, Y]_{\cal N_{\!{_1}}^{\boldsymbol{[}\stt\boldsymbol{]}}} \cap K(X) \\
=  & \hskip -5pt  A\Na(D[X], {\boldsymbol{[}\stt\boldsymbol{]}}) \cap K(X)
  \,.
\end{array}
$$
In particular, for all $E \in \FF(D)$,
$$
 \begin{array}{rl}
 E \Na(D, \star) \cap K = E^{\stt}=  & \hskip -5ptE^{\stt}[X] \cap K   \\
= & \hskip -5pt  (E[X])^{\boldsymbol{[}\stt\boldsymbol{]}}  \cap K \\
=  & \hskip -5pt
(E\Na(D[X], {\boldsymbol{[}\stt\boldsymbol{]}}) \cap K(X))   \cap K \\
=  & \hskip -5pt E\Na(D[X], {\boldsymbol{[}\stt\boldsymbol{]}})     \cap K\,.
\end{array}
$$

{\bf (2)}  Like ${\boldsymbol{[}\stt\boldsymbol{]}}$, also
 ${\boldsymbol{\langle}\stt\boldsymbol{\rangle}}$, and
 ${\boldsymbol{(}\stt\boldsymbol{)}}$  are finite type stable semistar operations on
 $D{_{\!{_1}}}:=D[X]$  (Propositions \ref{M}(2) and \ref{[star]}) then, for the following multiplicative subsets of $D{_{\!{_1}}}$,
 $$
 \begin{array}{rl}
  \cal S_{\!{_1}}^{\sharp}({\boldsymbol{\langle}\stt\boldsymbol{\rangle}})
  :=  & \hskip -6pt  \cal N_{\!{_1}}^{\boldsymbol{\langle}\stt\boldsymbol{\rangle}}  :=\{ 0 \neq  g{_{\!{_1}}} \in D{_{\!{_1}}}[Y] \mid \co_{D{_{\!{_1}}}} (g{_{\!{_1}}})^{\boldsymbol{\langle}\stt\boldsymbol{\rangle}} = D_{\!{_1}}^{\boldsymbol{\langle}\stt\boldsymbol{\rangle}} \}\,,\\
   \cal S_{\!{_1}}^{\sharp}({\boldsymbol{(}\stt\boldsymbol{)}}):= & \hskip -6pt  \cal N_{\!{_1}}^{\boldsymbol{(}\stt\boldsymbol{)}}
    := \{ 0\neq  g{_{\!{_1}}} \in D{_{\!{_1}}}[Y] \mid \co_{D{_{\!{_1}}}} (g{_{\!{_1}}})^{\boldsymbol{(}\stt\boldsymbol{)}}= D_{\!{_1}}^{\boldsymbol{(}\stt\boldsymbol{)}}  \}\,,
   \end{array}
   $$
   and, for all $A \in \FF(D[X])$, using also Proposition \ref{M}(1), we have:
    $$
 \begin{array}{rcl}
 A^{\boldsymbol{\langle}\stt\boldsymbol{\rangle}} = & \hskip -6pt
 A D[X, Y]_{\cal N_{\!{_1}}^{\boldsymbol{\langle}\stt\boldsymbol{\rangle}} } \cap K(X)  & \hskip -6pt =  A\Na(D[X], {\boldsymbol{\langle}\stt\boldsymbol{\rangle}}) \cap K(X)   \\
  & &  \hskip -6pt  =A\Na(D, \star) \cap AK[X] \,, \\
 A^{\boldsymbol{(}\stt\boldsymbol{)}} = & \hskip -6pt
 AD[X, Y]_{\cal N_{\!{_1}}^{\boldsymbol{(}\stt\boldsymbol{)}} } \cap K(X) & \hskip -6pt  = A \Na(D[X], {\boldsymbol{(}\stt\boldsymbol{)}}) \cap K(X)
 \\ & & \hskip -6pt  = A\Na(D, \star).
 \end{array}
 $$
 In particular, for all $E \in \FF(D)$, we have:
 $$
  \begin{array}{rl}
 E\Na(D, \star) \cap K = E^{\stt} = E^{\stt}[X] \cap K =
 & \hskip -6pt    (E[X])^{\boldsymbol{\langle}\stt\boldsymbol{\rangle}} \cap K =
 E\Na(D[X], {\boldsymbol{\langle}\stt\boldsymbol{\rangle}}) \cap K \\
  =& \hskip -6pt (E[X])^{{\boldsymbol{(}\stt\boldsymbol{)}}}   \cap K
   =  E\Na(D[X], {\boldsymbol{(}\stt\boldsymbol{)}} ) \cap K \,,
 \end{array}
 $$
 with $  (E[X])^{{\boldsymbol{(}\stt\boldsymbol{)}} } = E \Na(D, \star)$.

  \smallskip

Note that, by the previous descriptions of  ${\boldsymbol{\langle}\stt\boldsymbol{\rangle}}$ and ${\boldsymbol{(}\stt\boldsymbol{)}}$, we have:
$$
\begin{array}{rl}
P_1 \in \QSpec^{\boldsymbol{\langle}\stt\boldsymbol{\rangle}}(D[X]) & \Leftrightarrow \hskip 6pt P_1\Na(D, \star) \cap P_1K[X] \cap D[X] =\hskip 2pt  P_1\,,\\
P_1 \in \QSpec^{\boldsymbol{(}\stt\boldsymbol{)}}(D[X]) & \Leftrightarrow \hskip 6pt P_1\Na(D, \star)\cap D[X] =\hskip 2pt  P_1\,.
\end{array}
$$
Since $\Na(D, \star) = D[X]_{{\mathcal N}^\star}$, where ${\mathcal N}^\star:= \{ g \in D[X] \mid 0 \neq g \mbox{ and } \co_D(g)^\star = D^\star \}$, then:
 $$
\begin{array}{rl}
\QSpec^{\boldsymbol{\langle}\stt\boldsymbol{\rangle}}(D[X]) & \hskip -6 pt=\{ P_1 \in \Spec(D[X]) \mid \co_D(P_1)^{\stf} \subsetneq D^\star \mbox{ or } P_1\cap D =(0) \}\,, \\
\QSpec^{\boldsymbol{(}\stt\boldsymbol{)}}(D[X]) & \hskip -6 pt=\{ P_1 \in \Spec(D[X]) \mid \co_D(P_1)^{\stf} \subsetneq D^\star \}\,.
\end{array}
$$
Therefore, since   $\Max(\Na(D, \star)) = \{ Q[X] \mid Q \in \QMax^{\stf}(D)\}$ \cite[Proposition 3.1((2) and (3))]{FL}, we can conclude that
 $$
 \begin{array}{rl}
  \QMax^{\boldsymbol{\langle}\stt\boldsymbol{\rangle}} (D[X]) =&\hskip -7pt
 \{ Q_{\!{_1}} \in \Spec(D[X]) \mid  \!
Q_{\!{_1}} = Q[X] \mbox{ for some } Q \in \QMax^{\stf}(D) \} \  \cup \\
 & \hskip 8pt  \{(0)\neq Q_{\!{_1}}\in \Spec(D[X]) \mid  \!   Q_{\!{_1}}  \cap D = (0) \mbox{ and } \co_D(Q_{\!{_1}})^{\stf} = D^\star \}\,,\\
   \QMax^{\boldsymbol{(}\stt\boldsymbol{)}} (D[X]) =&\hskip -7pt
 \{ Q_{\!{_1}} \in \Spec(D[X]) \mid  \! Q_{\!{_1}}=
 Q[X] \mbox{ for some } Q \in \QMax^{\stf}(D) \}\,.
 \end{array}
 $$
\end{remark}

We already observed that the construction described in Remark \ref{sahandi}(1) is a modi\-fication of a previous construction due to Chang and Fontana \cite[Theorem 2.3]{CF}.
More precisely, Chang and Fontana considered the following
subset of $\Spec(D_1)$:
$$
\begin{array}{rl}
{\Delta}_{\!{_1}}^\prime := {\Delta}_{\!{_1, \star}}^\prime :=  \{Q_{\!{_1}} \in \Spec(D _{\!{_1}}) \mid &\mbox{either } \;\; Q_{\!{_1}} \cap D = (0) \;\;  \mbox{ or } \\
 & Q_{\!{_1}} =
(Q_{\!{_1}} \cap D)[X] \mbox{ and } (Q_{\!{_1}} \cap D)^{\stf} \subsetneq D^{\star}\}\,.
\end{array}
$$
Then, they considered the following associated extended saturated  multiplicative system in $D_{\!{_1}}[Y]$:
$$\mathcal{S}_{\!{_1}}^\prime:= \mathcal{S}_{\!{_1, \star}}^\prime := \mathcal{S}({\Delta}_{\!{_1, \star}}^\prime) :=  D_{\!{_1}}[Y] \setminus
\left(\bigcup\{Q_{\!{_1}}[Y] \mid Q_{\!{_1}} \in  {\Delta}_{\!{_1, \star}}^\prime \} \right)$$
 and the stable semistar operation of finite type of $D_{\!{_1}}$ defined by
$$
A^{\circlearrowleft_{\mathcal{S}_{\!{_{1, \star}}}^\prime }}:= AD_{\!{_1}} [Y]_{\mathcal{S}_{\!{_1}}^\prime} \cap K_{\!{_1}}\,, \;\;\;  \mbox{ for all } A \in
\FF(D_{\!{_1}}).
$$
They proved that, when $\star$ is  the $v$--operation (or the $t$--operation, or the $w$--operation)
 on $D$,
  then ${\circlearrowleft_{\mathcal{S}_{\!{_{1,\star}}}^\prime }}$ coincides with the $w$--operation of $D[X]$ \cite[Theorem 2.3(f)]{CF}. (Note that, in that paper, the authors denoted the semistar operation ${\circlearrowleft_{\mathcal{S}_{\!{_{1,\star}}}^\prime }}$ of $D[X]$ by $[\star]$; we avoid now this notation, since we already use it here with a different meaning.)

\begin{proposition} \label{<star>}
Let $\star$ be a semistar operation   defined on an integral domain $D$.  Let $\boldsymbol{\star}_{\!{_1}}^\prime := \circlearrowleft_{\mathcal{S}_{\!{_{1, \star}}}^\prime }$ be the stable semistar operation of finite type on $D[X]$ defined above and let
$ \boldsymbol{\star}_{\!{_1}} $ be  the semistar operation  defined in Remark \ref{sahandi}(1)
(i.e., $ \boldsymbol{\star}_{\!{_1}}   := {\circlearrowleft_{\mathcal S^{\sharp}_{\!{_{1, \star}}}}}\!  =  \boldsymbol{[}\stt\boldsymbol{]}$).
\begin{enumerate}
\item[{\bf (1)}]
  $ \boldsymbol{[}\stt\boldsymbol{]} =\boldsymbol{\star}_{\!{_1}} \leq \boldsymbol{\star}_{\!{_1}}^\prime$.
 \item[{\bf (2)}]
 $\boldsymbol{\star}_{\!{_1}}^\prime = {\boldsymbol{\langle}\stt\boldsymbol{\rangle}}$.
 \end{enumerate}
\end{proposition}
\begin{proof}
(1) follows easily from the fact that  $ {\Delta}_{\!{_1, \star}} \supseteq {\Delta}_{\!{_1, \star}}^\prime$
  or, equivalently,
  ${\mathcal S}^{\sharp}_{\!{_1}} \subseteq {\mathcal{S}}_{\!{_1}}^{\prime} $.

  (2) We know that  $\QMax^{\boldsymbol{\star}_{\!{_1}}^\prime}(D_{\!{_1}}) = \{Q_{\!{_1}}   \in \Spec(D_{\!{_1}}) \mid  Q_{\!{_1}} \cap D = (0)$ and
$\co_D(Q_{\!{_1}})^{\stf} = D^\star \} \cup \{Q[X] \mid  Q \in  \QMax^{\stf}(D)\}$ \cite[Theorem 2.3(e)]{CF}. On the other hand, by what we observed in Remark \ref{sahandi}(2),
$$
\QMax^{\boldsymbol{\star}_{\!{_1}}^\prime}(D[X]) = \QMax^{\boldsymbol{\langle}\stt\boldsymbol{\rangle}} (D[X])\,,$$
and so, since   $\boldsymbol{\star}_{\!{_1}}^\prime $ and $ {\boldsymbol{\langle}\stt\boldsymbol{\rangle}}$ are both stable semistar operations of finite type, we conclude that $ \boldsymbol{\star}_{\!{_1}}^\prime = {\boldsymbol{\langle}\stt\boldsymbol{\rangle}}$.
  \end{proof}

  As a final remark (with the notation of the present paper),  note that in \cite[Corollary 2.5(1)]{CF} the authors prove that $D$ is a Pr\"ufer $\star$--multiplication domain if and only if $D[X]$ is a Pr\"ufer ${\boldsymbol{\langle}\stt\boldsymbol{\rangle}}$--multiplication domain. On the other hand, it is not true that
$D$ is a Pr\"ufer $\star$--multiplication domain if and only if $D[X]$ is a Pr\"ufer ${\boldsymbol{[}\stt\boldsymbol{]}}$--multiplication domain (take, for instance, $D$ a Pr\"ufer domain, but not a field, and $\star =d_D$, in this case
${\boldsymbol{[}\widetilde{\ d_D}\boldsymbol{]}} = {\boldsymbol{[}d_D\boldsymbol{]}} = d_{D[X]}$ and, obviously, $D[X]$ is not a Pr\"ufer domain).
This fact justifies the terminology used in \cite{CF}, where the authors call the semistar operation denoted here by ${\boldsymbol{\langle}\stt\boldsymbol{\rangle}}$ the stable semistar operation of finite type canonically associated to $\star$.

\medskip
\noindent
{\bf Acknowledgment.} We thank Muhammad Zafrullah for the helpful suggestions received during the preparation of the present paper and the referee  for providing constructive comments. \\
The first author's work was supported by the 2008 Research Fund from the College of Natural
Sciences, University of Incheon.
\\The second author was partially supported by a MIUR-PRIN grant N. 2008WYH9NY.

\medskip
\end{document}